\documentclass{amsart}

\usepackage{amsmath, amssymb, amsthm}
\usepackage{graphicx}
\usepackage{graphics}
\usepackage{epsfig}
\usepackage{color}
\usepackage[all]{xy}
\usepackage{amsfonts}
\usepackage{latexsym}
\usepackage{hyperref}

\usepackage[normalem]{ulem} 
\usepackage{color}

\definecolor{darkgreen}{rgb}{0,0.5,0}

\newtheorem{defin}{Definition}[section]
\newtheorem{thm}[defin]{Theorem}
\newtheorem{cor}[defin]{Corollary}
\newtheorem{rem}[defin]{Remark}
\newtheorem{lem}[defin]{Lemma}
\newtheorem{prop}[defin]{Proposition}
\newtheorem{quest}[defin]{Question}
\newtheorem{prob}[defin]{Problem}
\newtheorem{axi}[defin]{Axiom}
\newtheorem{ex}[defin]{Example}

\newtheorem{conven}[defin]{Convention}
\newtheorem{conj}[defin]{Conjecture}

\newcommand{\theorem}[1]{\begin{thm} #1 \end{thm}}
\newcommand{\theoremname}[2]{\begin{thm}[#1] #2 \end{thm}}
\newcommand{\proposition}[1]{\begin{prop} #1 \end{prop}}
\newcommand{\propositionname}[2]{\begin{prop}[#1] #2 \end{prop}}
\newcommand{\lemma}[1]{\begin{lem} #1 \end{lem}}
\newcommand{\definition}[1]{\begin{defin} \emph{#1} \end{defin}}
\newcommand{\question}[1]{\begin{quest} \emph{#1} \end{quest}}
\newcommand{\remark}[1]{\begin{rem} \emph{#1} \end{rem}}
\newcommand{\remarkname}[2]{\begin{rem}[#1] \emph{#2} \end{rem}}
\newcommand{\corollary}[1]{\begin{cor} #1 \end{cor}}

\newcommand{\lemmaname}[2]{\begin{lem}[#1] #2 \end{lem}}

\newcommand{\convention}[1]{\begin{conven} \emph{#1} \end{conven}}
\newcommand{\conjecture}[1]{\begin{conj} #1 \end{conj}}
\newcommand{\ov}[1]{\overline{#1}}
\newcommand{\wh}[1]{\widehat{#1}}
\newcommand{\wt}[1]{\widetilde{#1}}
\newcommand{\PL}{\mathrm{PL}}
\newcommand{\EP}{\mathrm{EP}}
\newcommand{\REP}{\mathcal{R}\cdot \EP_2}
\newcommand{\WTEP}{\widetilde{\mathrm{EP}}}
\newcommand{\id}{\mathrm{id}}
\newcommand{\Aut}{\mathrm{Aut}}
\newcommand{\Inn}{\mathrm{Inn}}
\newcommand{\Fix}{\mathrm{Fix}}
\newcommand{\R}{\mathbb{R}}

\newcommand{\supp}{\operatorname{supp}}

\DeclareMathSymbol{\Rb}{\mathbin}{AMSb}{"52}
\DeclareMathSymbol{\Zb}{\mathbin}{AMSb}{"5A}

\title{\textsc{The conjugacy problem in extensions of Thompson's group $F$}}

\author{Jos\'e Burillo}
\address{Departament de Matem\`atica Aplicada IV, Escola Polit\`ecnica
Superior de Castelldefels, Universitat Polit\`ecnica de Catalunya, C/Esteve
Torrades 5, 08860 Castelldefels, Barcelona, Spain}
\email{burillo@ma4.upc.edu}
\author{Francesco Matucci}
\address{D\'epartement de Math\'ematiques, Facult\'e des Sciences d'Orsay,
Universit\'e Paris-Sud 11, B\^atiment 425, Orsay, France}
\email{francesco.matucci@math.u-psud.fr}
\author{Enric Ventura}
\address{Departament Matem\`atica Aplicada III, Universitat Polit\`ecnica de Catalunya, Manresa,
Catalunya}
\email{enric.ventura@upc.edu}

\thanks{The first and third author acknowledge support from the MEC grant MTM2011-25955.
The second author gratefully acknowledges the Fondation Math\'ematique Jacques Hadamard (FMJH - ANR - Investissement d'Avenir)
and its staff for the support received during the development of this work.
The three authors gratefully
acknowledge the Centre de Recerca Matem\`atica (CRM) and its staff for the support received during the development of
this work.}

\begin{document}

\maketitle

\begin{abstract}
We solve the twisted conjugacy problem on Thompson's group $F$. We also exhibit orbit
undecidable subgroups of $\Aut(F)$, and give a proof that $\Aut(F)$ and $\Aut_+(F)$ are
orbit decidable provided a certain conjecture on Thompson's group $T$ is true. By using
general criteria introduced by Bogopolski, Martino and Ventura in~\cite{bomave2}, we
construct a family of free extensions of $F$ where the conjugacy problem is unsolvable.
As a byproduct of our techniques, we give a new proof of a result of
Bleak-Fel'shtyn-Gon\c{c}alves in~\cite{bleakfelshgonc1} showing that $F$ has property
$R_\infty$, and which can be extended to show that Thompson's group $T$ also has property
$R_\infty$.
\end{abstract}


\section{Introduction \label{sec:intro}}

Since Max Dehn formulated the three main problems in group theory in 1911, they have been a
central subject of study in the theory of infinite groups. 
There now exists a large body of works devoted to the study of these problems.
In this paper we focus on the conjugacy problem and a variant
known as \emph{the twisted conjugacy problem}.
The conjugacy problem is known to be solvable for Thompson's groups
$F,T$ and $V$ by works of Guba and Sapir \cite{gusa1}, Belk and the second author \cite{matucci9}
and Higman \cite{hig}.
Our interest arose in the study of extensions of the group $F$ where we find an unsolvability result.
Even though Thompson himself used the groups $F,T,V$ in the construction of finitely presented
groups with unsolvable word problem, to the best of our knowledge, 
the result that we obtain is a first in a direct generalization of the original Thompson groups.
Moreover, we also look at property $R_\infty$ which has been under study recently
and which is known to true for the group $F$ and one of its extensions.

We now give a more detailed description of the results.
Let $F$ be a group. We say that a subgroup $A\leqslant \Aut(F)$ has solvable \emph{orbit
decidability problem} (ODP) if it is decidable to determine, given $y,z \in F$, whether or
not there is $\varphi \in A$ and $g \in F$ such that
 $$
\varphi(z)= g^{-1} y g.
 $$

On the other hand, if $\varphi \in \Aut(F)$, we say that $F$ has solvable
\emph{$\varphi$-twisted conjugacy problem} (TCP$_\varphi$) if it is decidable to determine,
given $y,z \in F$, whether or not they are $\varphi$-\emph{twisted conjugated} to each
other, i.e. whether there exists $g\in F$ such that
\begin{equation}\label{eq:TCP-equation}
\numberwithin{equation}{section} z = g^{-1} y \varphi(g).
\end{equation}
More generally, we say that the group $F$ has solvable \emph{twisted conjugacy problem}
(TCP) if (TCP$_\varphi$) is solvable for any given $\varphi \in \Aut(F)$.

In their recent paper~\cite{bomave2}, Bogopolski, Martino and Ventura developed a criterion
to study the conjugacy problem for some extensions of groups, and found a connection of
this problem with the two problems mentioned above.

Let $F,G,H$ be finitely presented groups and consider a short exact sequence
\begin{equation}\label{eq:usual-exact-sequence}
\numberwithin{equation}{section} 1 \longrightarrow F \overset{\alpha}{\longrightarrow} G
\overset{\beta}{\longrightarrow} H \longrightarrow 1.
\end{equation}
In this situation, $\alpha(F) \unlhd G$ and so the conjugation map $\varphi_g$, for $g \in
G$, restricts to an automorphism of $F$, $\varphi_g \colon F\to F$, $x\mapsto g^{-1}xg$,
(which does not necessarily belong to $\Inn(F)$). We define the \emph{action subgroup} of
the sequence (\ref{eq:usual-exact-sequence}) to be the group of automorphisms
 $$
A_G =\{\varphi_g \mid g \in G\} \leqslant \Aut(F).
 $$

\theoremname{Bogopolski-Martino-Ventura, \cite{bomave2}}{ \label{thm:bomave-extensions}
Let
 $$
1 \longrightarrow F \overset{\alpha}{\longrightarrow} G \overset{\beta}{\longrightarrow} H \longrightarrow 1.
 $$
be an algorithmic short exact sequence of groups such that
\begin{enumerate}
\item $F$ has solvable twisted conjugacy problem,
\item $H$ has solvable conjugacy problem, and
\item for every $1 \ne h \in H$, the subgroup $\langle h \rangle$ has finite index in its
    centralizer $C_H(h)$, and we can compute a set of coset representatives of $\langle
    h\rangle$ in $C_H(h)$.
\end{enumerate}
Then, the conjugacy problem for $G$ is solvable if and only if the action subgroup $A_G
=\{\varphi_g \mid g \in G\}\leqslant \Aut(F)$ is orbit decidable.
}

Here, a short exact sequence is \emph{algorithmic} if all the involved groups are finitely
presented and given to us with an explicit finite presentation, and all the morphisms are
given by the explicit images of the generators.

Condition (3) is of more technical nature. It is clearly satisfied in free groups (where
the centralizer of a non-trivial element $h$ is just the cyclic subgroup generated by its
maximal root $\hat{h}$), and it is also true in torsion-free hyperbolic groups,
see~\cite{bomave2}.

The goal of the present paper is to study the conjugacy problem in some extensions of
Thompson's group $F$ via Theorem~\ref{thm:bomave-extensions} (see~\cite{bomave2, ventura1} for
references to similar applications of this same theorem into other families of groups).

\bigskip

We will assume the reader is familiar with Thompson's groups $F$ (also denoted by
$\PL_2(I)$, where $I=[0,1]$ is the unit interval) and $T$ (also denoted by $\PL_2(S^1)$,
where $S^1$ is the unit circle) and in any case, the comprehensive survey by Cannon, Floyd
and Parry \cite{cfp} is an excellent source of information for Thompson's groups.

We will employ techniques on conjugacy in the Bieri-Thompson-Stein-Strebel groups used by
Kassabov and the second author in \cite{matucci5} and a rephrasing by Belk and the second
author in \cite{matucci9, matucci3} of a conjugacy invariant of Brin and Squier
\cite{brin2}. The idea is to assume that the twisted conjugacy equation has a solution and
use this to determine necessary conditions that a twisted conjugator should satisfy. This
allows one to build some candidate conjugators which must then be tested.

With these techniques, we obtain the first result in the paper:

\theorem{\label{thm:TCP-solvable}
Thompson's group $F$ has solvable twisted conjugacy problem.
}

Putting together Theorems~\ref{thm:bomave-extensions} and~\ref{thm:TCP-solvable}, this
opens us to the possibility of finding extensions of $F$ with solvable/unsolvable conjugacy
problem, by detecting subgroups of $\Aut(F)$ which are orbit decidable/orbit undecidable:

\theorem{\label{coro}
Consider Thompson's group $F=\PL_2(I)$, a torsion-free hyperbolic group $H$, and let
 $$
1 \longrightarrow F \overset{\alpha}{\longrightarrow} G \overset{\beta}{\longrightarrow} H \longrightarrow 1.
 $$
be an algorithmic short exact sequence. The group $G$ has solvable conjugacy problem if and
only if the action subgroup $A_G \leqslant \Aut (F)$ is orbit decidable.
}

Using the previous result one can create extensions of $F$ with unsolvable conjugacy
problem.

\theorem{\label{thm:CP-extension-unsolvable}
There are extensions of Thompson's group $F$ by finitely generated free groups, with
unsolvable conjugacy problem.
} It is also possible to build some interesting extensions of $F$ with solvable conjugacy
problem, provided that an open conjecture about $F$ is true. We study this in Section
\ref{sec:ODP-solvable}.

\bigskip

A group $G$ has the \emph{property $R_\infty$} if it has infinitely many distinct
$\varphi$-twisted conjugacy classes, for any $\varphi \in \Aut(G)$. Thompson's group $F$
was shown to have property $R_\infty$ by Bleak, Fel'shtyn and Gon\c{c}alves in
\cite{bleakfelshgonc1}. We give an alternative proof, which can be extended to Thompson's
group $T$.

\theorem{\label{thm:R-infty-T}
Thompson's group $T$ has property $R_\infty$.
}

The paper is organized as follows. In Section~\ref{sec:twisted-problem} we introduce the
groups we will be working with, we restate the twisted conjugacy problem for $F$ and prove
Theorems~\ref{thm:TCP-solvable} and~\ref{coro}. In Section~\ref{sec:CP-extensions} we
construct orbit undecidable subgroups of $\Aut(F)$ and exhibit free extensions of $F$ with
unsolvable conjugacy problem. In Section~\ref{sec:ODP-solvable} we consider orbit
decidability and construct some interesting extensions of $F$, which happen to have
solvable conjugacy problem assuming an open conjecture on $F$ is true. In
Section~\ref{ssec:R-infty} we show that the groups $F$ and $T$ have property $R_{\infty}$
using ideas from Section~\ref{sec:twisted-problem}. Finally, in
Section~\ref{sec:generaltions-of-results} we analyze the extent to which the techniques of
this paper generalize to other families of Thompson-like groups.

\subsection*{Acknowledgments}

The authors would like to thank Matt Brin, Collin Bleak, Martin Kassabov, Jennifer Taback
and Nathan Barker for helpful conversations about this work.

\section{\label{sec:twisted-problem} The twisted conjugacy problem for $F$}

In this section we prove Theorem~\ref{thm:TCP-solvable}. The techniques developed for this
purpose will be later used in Section~\ref{ssec:R-infty} to obtain a couple of byproducts.

\subsection{Thompson's group and its automorphisms\label{sec:thompson and autos}}

We will look at Thompson's group $F$ from different perspectives. The standard one is to
look at $F$ as the group $\PL_2(I)$ of orientation preserving piecewise-linear
homeomorphisms of the unit interval $I=[0,1]$ with a discrete (and hence finite) set of
breakpoints at dyadic rational points, and such that all slopes are powers of $2$ (the
interval $I$ can be replaced to an arbitrary $[p,q]$ with $p,q$ being dyadic rationals and
the resulting group is clearly isomorphic). We will also need to regard $F$ as a subgroup
of a bigger group: consider the group $\PL_2(\mathbb{R})$ of all piecewise-linear
homeomorphisms of $\mathbb{R}$ with a discrete set of breakpoints at dyadic rational points
and such that all slopes are powers of $2$; and consider the subgroup of those elements $f$
which are eventually integral translations, i.e. for which there exist $m_-, m_+\in
\mathbb{Z}$ and $L,R\in \mathbb{R}$ such that $f(x)=x+m_-$ for all $x\leqslant L$, and
$f(x)=x+m_+$ for all $x\geqslant R$. It is straightforward to see that this subgroup of
$\PL_2(\mathbb{R})$ is isomorphic to $\PL_2(I)$; see Proposition~3.1.1 in Belk and
Brown~\cite{bebr} for an explicit isomorphism (it is interesting to note that, through this
isomorphism, $2^{m_-}$ is the slope at the right of 0, and $2^{m_+}$ the slope at the left
of 1). Both copies of Thompson's group will be denoted $F$, and it will be clear from the
context which one are we talking about at any moment.

Thompson's group admits a finite presentation. The two generators are usually written $x_0$
and $x_1$, which represent the following maps on the real line:
$$
x_0(t)=t+1\qquad
x_1(t)=\left\{\begin{array}{ll}
t&\text{if }t<0\\
2t&\text{if }0\leq t\leq 1\\
t+1&\text{if }t>1.
\end{array}\right.
$$
With these generators, $F$ admits a finite presentation with just two relations, which have
lengths 10 and 14. See \cite{cfp} for details. Moreover, as we will need this later, we observe that
when we regard $F$ as the group $\PL_2([0,1])$, the generator $x_0$ has this form:
\[
\theta(t):=
\begin{cases}
2t & t \in \left[0,\frac{1}{4} \right] \\
t + \frac{1}{4}& t \in \left[\frac{1}{4},\frac{1}{2} \right] \\
\frac{t}{2}+\frac{1}{2} & t \in \left[\frac{1}{2},1 \right].
\end{cases}
\]
We distinguish $x_0$ and $\theta$ to make it clear that the first one is seen as an element of $\PL_2(\mathbb{R})$
while the second is regarded as a map in $\PL_2([0,1])$.
The \emph{support} of an element $f\in \PL_2(\mathbb{R})$ is the collection of points where
it is different from the identity, $\supp (f)=\{ t\in \mathbb{R}\mid f(t)\neq t\}$.

\definition{\label{thm:definition-subgroups}
We define the following subgroups of $\PL_2(\mathbb{R})$:
\begin{enumerate}
\item $\EP_2=\{f \in \PL_2(\mathbb{R}) \mid \exists\, L,R \in \mathbb{R} \text{ such that
    } f(t-1)=f(t)-1 \,\, \forall t\leqslant L, \text{ and } f(t+1)=f(t)+1\,\, \forall
    t\geqslant R \}$ i.e., all functions in $\PL_2(\mathbb{R})$ which are ``eventually
    periodic'' and orientation preserving.
\item $F=\{f \in \EP_2 \mid \exists\, L,R \in \mathbb{R},\, \exists\, m_-,m_+ \in
    \mathbb{Z} \text{ such that } f(t)=t+m_-\,\, \forall t\leqslant L, \text{ and }
    f(t)=t+ m_+\,\, \forall t\geqslant R\}$. As noted above, $F\simeq \PL_2(I)$ is the
    standard copy of Thompson's group inside $\PL_2(\mathbb{R})$.
\item Let $G$ be any subset of $\EP_2$. For every $-\infty\leqslant p<q\leqslant
    +\infty$, define $G(p,q)$ to be the set of elements in $G$ with support inside the
    interval $(p,q)$ i.e. $G(p,q)=\{ g\in G\mid g(t)=t, \forall t \not \in (p,q)\}$ (so,
    $G(-\infty, +\infty)=G$). Also, define $G^>=\{g \in G \mid g(t)>t, \forall t \in
    \mathbb{R}\}$ and, similarly, $G^<$. When combining both notations we shall
    understand the inequality restricted to the support, i.e. $G^<(p,q)=\{ g\in G\mid
    g(t)=t,\,\, \forall t \not \in (p,q), \text{ and } g(t)<t\,\, \forall t\in (p,q)\}$.
    Note that, if $G$ is a subgroup, then $g\in G^>(p,q)$ if and only if $g^{-1}\in
    G^<(p,q)$.
\end{enumerate}
At certain point in the arguments we will also need to consider orientation reversing maps.
Admitting both orientations in the definition above, one can define the group
$\PL_2^{\pm}(\mathbb{R})$ and the corresponding subgroup $\WTEP_2 =\{f \in
\PL_2^{\pm}(\mathbb{R}) \mid \exists\, L,R \in \mathbb{R}, \exists\, \epsilon =\pm 1 \text{
such that } f(t-1) =f(t)-\epsilon \,\, \forall t\leqslant L, \text{ and } f(t+1)=f(t)
+\epsilon\,\, \forall t\geqslant R \}$. Note that $\EP_2$ is a subgroup of $\WTEP_2$ of
index two, and $\WTEP_2 =\EP_2 \cup \mathcal{R}\cdot \EP_2$, where $\mathcal{R} \in \WTEP_2
\setminus \EP_2$ is the \emph{reversing map}, $\mathcal{R}(t)=-t$ for all $t\in
\mathbb{R}$.}

\begin{center}
\begin{minipage}{0.4\columnwidth}
\centering
\includegraphics[height=5.5cm]{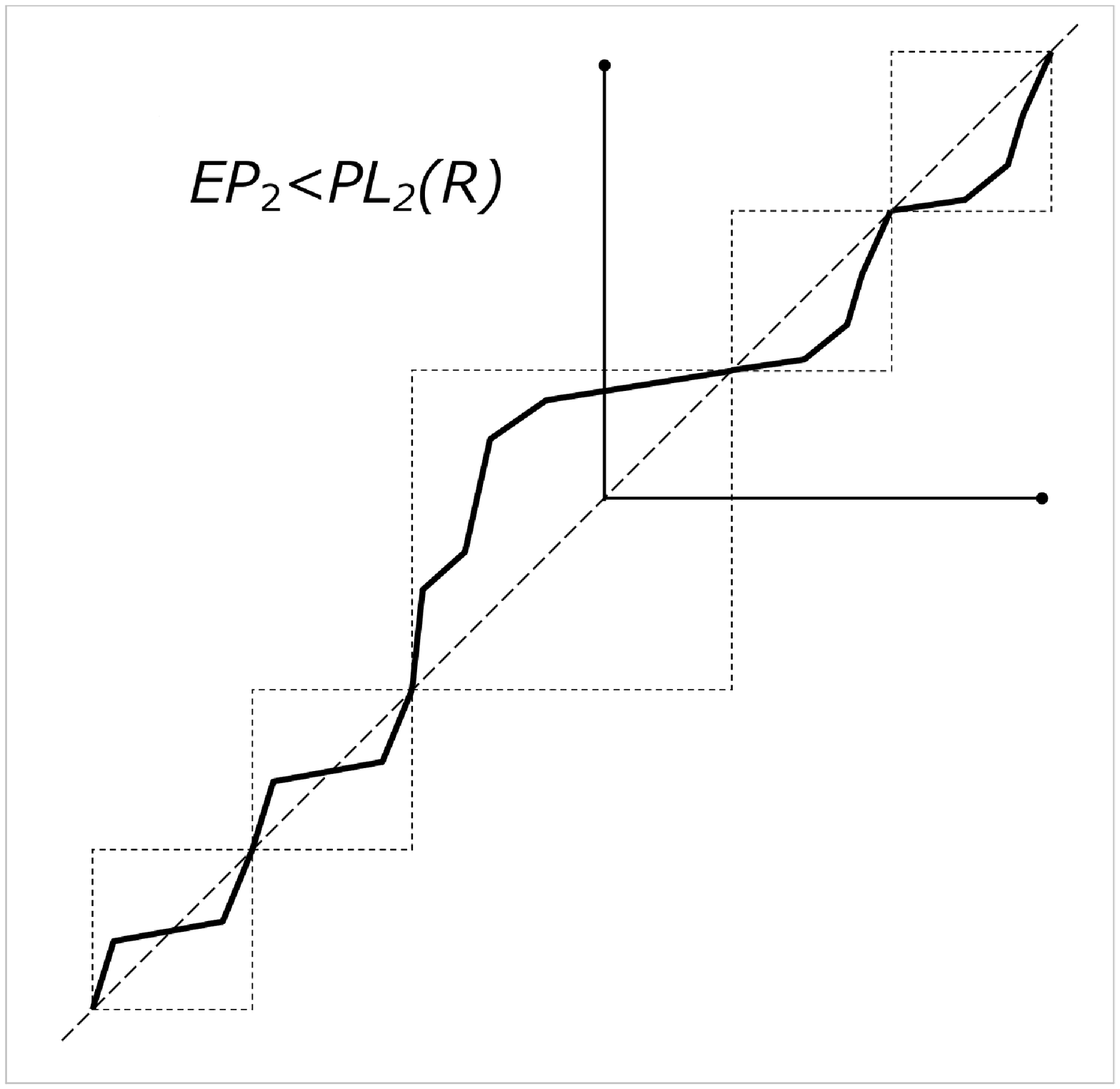}\\
\end{minipage}%
\hspace{0.3in}
\begin{minipage}{0.4\columnwidth}
\centering
\includegraphics[height=5.5cm]{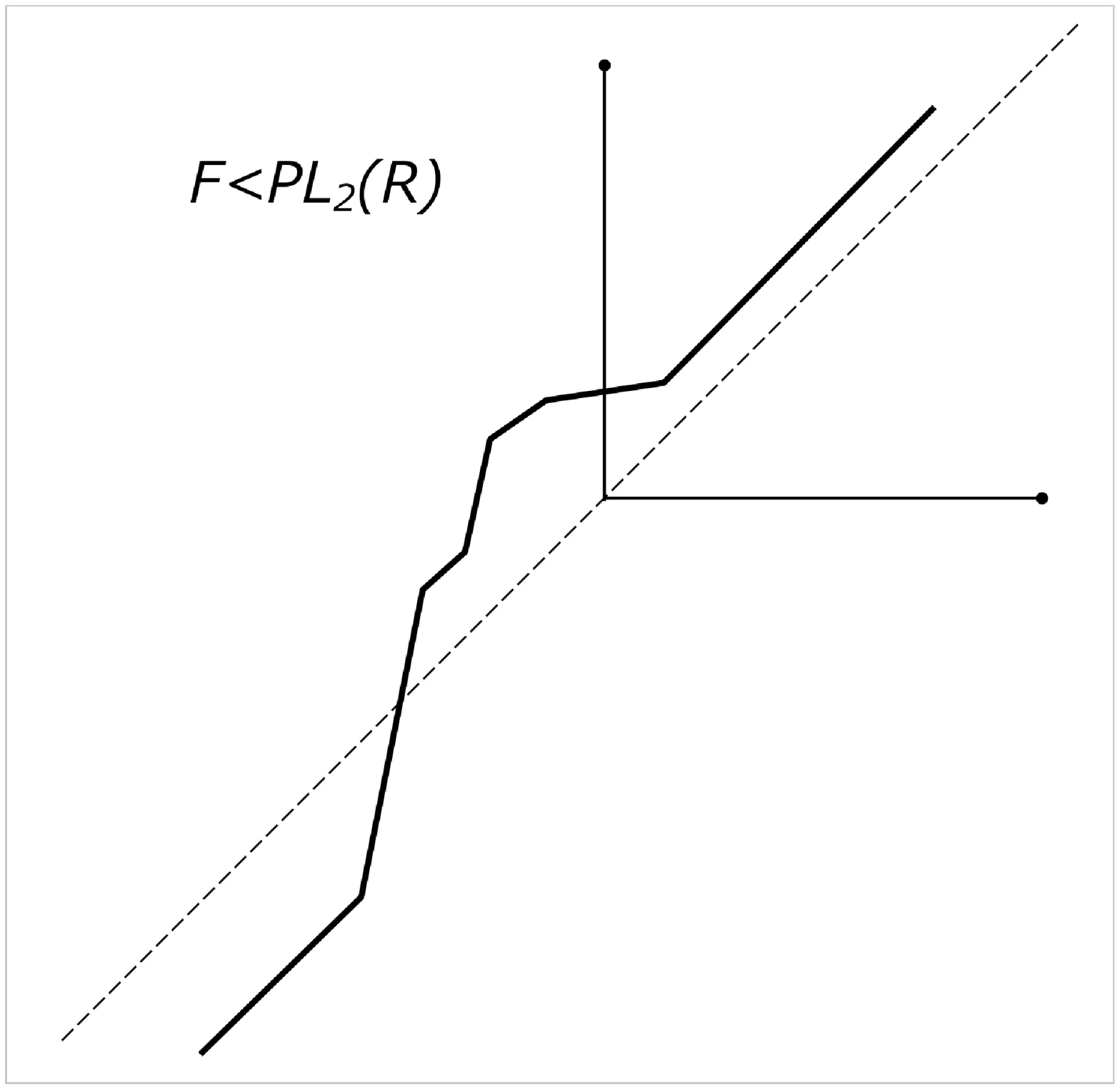}\\
\end{minipage}%
\end{center}

\convention{
When talking about elements $f\in \PL_2 (\R)$, we say that a property $\mathcal{P}$ holds
\emph{for $t$ positive sufficiently large} (respectively, \emph{for $t$ negative
sufficiently large}) to mean that there exists a number $R>0$ such that $\mathcal{P}$ holds
for every $t\geqslant R$ (respectively, there exists a number $L<0$ such that $\mathcal{P}$
holds for every $t\leqslant L$). For example, $f\in F$ if and only if it is an integral
translation for $t$ positive sufficiently large, and for $t$ negative sufficiently large.
}

\remark{\label{thm:remark-initial-slope}
Observe that, for $g\in F\leqslant \PL_2(\mathbb{R})$, the integer $m_-$ above (satisfying
that $g(t)=t+m_-$ for $t$ negative sufficiently large) can also be obtained as the limit
$m_- =\lim_{t\to -\infty} g(t)-t$. Similarly, $g(t)=t+m_+$ for $t$ positive sufficiently
large, where $m_+ =\lim_{t\to +\infty} g(t)-t$. These two real numbers are called,
respectively, the \emph{initial slope} and the \emph{final slope} of $g$ because, when
regarded as an element of $\PL_2(I)$, the slopes of $g$ on the right of the point 0 and on
the left of the point 1 are, precisely, $2^{m_-}$ and $2^{m_+}$, respectively.
}

\subsection{Automorphisms and transitivity on dyadics}

To deal with the $\varphi$-twisted conjugacy problem for $F$, we first need to understand
what the automorphisms of Thompson's group $F$ look like. They have all been classified by
Brin in his Theorem~1 in~\cite{brin5} (see also Theorem~1.2
in~\cite{burcleary2} for a more explicit version). The key idea to
understand $\Aut (F)$ is the fact that conjugation by elements from $\WTEP_2$ preserves
$F$, and these conjugations give precisely all automorphisms of $F$:


\theoremname{Brin, \cite{brin5}}{\label{thm:brin-thm}
For Thompson's group $F$, the map
 $$
\begin{array}{ccc}
\WTEP_2 & \longrightarrow & \Aut (F) \\ \tau & \mapsto & \begin{array}[t]{rcl} \gamma_{\tau}\colon F & \rightarrow & F
\\ g & \mapsto & \tau^{-1}g\tau, \end{array}
\end{array}
 $$
is well defined and it is a group isomorphism, so $\Aut (F) \simeq \WTEP_2$. Furthermore,
given $\varphi \in \Aut (F)$ by the images of the standard generators, one can
algorithmically compute the (unique) $\tau\in \WTEP_2$ such that $\varphi
(g)=\tau^{-1}g\tau$ for all $g\in F$.
}

\definition{\label{aut+}
We denote by $\Aut_+ (F)$ the group of automorphisms of $F$ given by conjugation by
orientation preserving $\tau$'s (see Theorem~\ref{thm:brin-thm}); it is an index two
subgroup $\EP_2 \simeq \Aut_+(F) <_2 \Aut (F)\simeq \WTEP_2$.
}


\remarkname{Explicit rewriting of elements of $\Aut(F)$}{ Theorem~\ref{thm:brin-thm},
including its algorithmic contents, is crucial for the arguments of the present paper.
Brin's original theorem establishes the isomorphism and we can do the algorithmic
determination of $\tau$ in the following form. Burillo and Cleary \cite{burcleary2}
obtained a finite presentation for $\Aut(F)$ with nine generators $\varphi_1\,\ldots ,
\varphi_9$ all expressed in terms of the standard presentation of $F$, and as conjugations
by suitable $\tau_1,\ldots ,\tau_9\in \WTEP_2$, i.e. $\varphi_i=\gamma_{\tau_i}$ for
$i=1,\ldots ,9$. Suppose $\varphi \in \Aut(F)$ is given by the images of $x_0$ and $x_1$.
We can enumerate all formal words $w$ on letters $\varphi_1,\ldots ,\varphi_9$ and for each
one compute the images of $x_0$ and $x_1$ by $w(\varphi_1,\ldots ,\varphi_9)$ until they
match with $\varphi(x_0)$ and $\varphi(x_1)$ (here we need to use the word problem for
$F$); this match will happen sooner or later because $\varphi_1,\ldots ,\varphi_9$ do
generate $\Aut(F)$. Once we have this word, it is clear that $\tau=w(\tau_1,\ldots
,\tau_9)\in \WTEP_2$ satisfies $\gamma_{\tau}=\gamma_{w(\tau_1,\ldots
,\tau_9)}=w(\gamma_{\tau_1},\ldots ,\gamma_{\tau_9})=w(\varphi_1,\ldots
,\varphi_9)=\varphi$.}

The following is a result explaining how to build $\PL_2$-maps acting in a prescribed way
on some given rational numbers. The first part gives an arithmetic condition for the
existence of such a map. The second part expresses the flexibility of these groups: one can
always ``cut" the graphical representation of an element at a given dyadic rational, and
freely ``glue" the pieces to obtain new elements. This result will often be needed along
the present paper.

\propositionname{Kassabov-Matucci, \cite{matucci5}}{ \label{thm:rationals-coincide} Let
$\eta,\zeta$ be dyadic rationals, let $\alpha, \beta\in \mathbb{Q}\cap(\eta,\zeta)$ written
in the form $\alpha=\frac{2^t m}{n}$ and $\beta=\frac{2^k p}{q}$ with $t,k\in \mathbb{Z}$
and $m,n,p,q$ odd integers such that $(m,n)=(p,q)=1$, and let $\eta <\alpha_1 <\cdots
<\alpha_r <\zeta$ and $\eta <\beta_1 <\cdots <\beta_r <\zeta$ be two finite sequences of
rational numbers.
\begin{enumerate}
\item The following are equivalent:
\begin{itemize}
\item[(a)] there exists $g\in \PL_2([\eta,\zeta ])$ such that $g(\alpha)=\beta$,
\item[(b)] there exists $g\in \PL_2 (\mathbb{R})$ such that $g(\alpha)=\beta$,
\item[(c)] there exists $g\in \EP_2$ such that $g(\alpha)=\beta$,
\item[(d)] there exists $g\in F$ such that $g(\alpha)=\beta$,
\item[(e)] $n=q$ and $p\equiv 2^Rm \pmod{n}$ for some $R\in \mathbb{Z}$.
\end{itemize}
Moreover, there is an algorithm which constructs such elements $g$ if condition (e) is
satisfied.
\item There exists $g\in F$ with $g(\alpha_i )=\beta_i$ if and only if for every
    $i=1,\ldots, r$ there exists $g_i \in F$ such that $g_i(\alpha_i) = \beta_i$.
    Moreover, if such a $g$ exists it can be constructed from the $g_i$'s.
\end{enumerate}
}

The following is a well known standard result (see for example \cite{matucci5} for a
proof).

\lemma{\label{thm:standard-folklore} Let $p \in \mathbb{Q}$ and $g \in \PL_2([p,p+1])$. Let
$u,v \in (p,p+1)$ be such that $u \not \in \Fix(g)$. Then there exists at most a unique
integer $m$ such that $g^m(u)=v$, and one can algorithmically decide it (and compute such
an $m$ if it exists).}

\subsection{Restatement of the TCP \label{sec:restatement-TCP}}

Our goal in this section is to solve the twisted conjugacy problem in $F$: given $\varphi
\in \Aut(F)$ and $y,z\in F$ (all in terms of the standard presentation of $F$, i.e.
$\varphi (x_0),\, \varphi (x_1),\, y,\, z$ are given to us as words on $x_0,\, x_1$), we
have to decide whether there exists $g\in F$ such that
 $$
z=g^{-1} y \varphi(g).
 $$

Applying Theorem~\ref{thm:brin-thm}, we can compute $\tau \in \WTEP_2$ such that $\varphi
(g)=\tau^{-1}g\tau$ for all $g\in F$, and the previous equation becomes $z=g^{-1} y
(\tau^{-1}g\tau)$, that is
 $$
z \tau^{-1} = g^{-1} (y \tau^{-1})g.
 $$
Relabeling $\ov{y}:=y\tau^{-1}\in \WTEP_2$ and $\ov{z}:= z\tau^{-1}\in \WTEP_2$ to get
\begin{equation}\label{eq:conjugacy}\numberwithin{equation}{section}
\ov{z} = g^{-1} \ov{y} g,
\end{equation}
the problem reduces to the standard conjugacy problem in $\WTEP_2$, but with the conjugator
$g$ forced to be chosen from $F\leqslant \WTEP_2$.

\definition{
Given two elements $\ov{y}, \ov{z}\in \WTEP_2$, we write $\ov{y}\sim_F \ov{z}$ if they are
conjugated by a conjugator in $F$, i.e. if there exists $g\in F$ such that $\ov{z} = g^{-1}
\ov{y} g$.
}

Notice that if one of $\ov{y}$ and $\ov{z}$ is in $\EP_2$ and the other is not, then
equation~(\ref{eq:conjugacy}) has no solution. Thus, we can split its study into two cases:
the orientation preserving case, i.e. when $\ov{y},\ov{z} \in \EP_2$ (studied in
Sections~\ref{ssec:periodicity-boxes}, \ref{ssec:fixed-points},
\ref{ssec:reducing-to-squares} and \ref{sec:rescaling-the-circle}) and then the orientation
reversing one, i.e. when $\ov{y},\ov{z}\in \REP$ (considered in
Section~\ref{sec:special-case}). Finally, in Section~\ref{ssec:solution-TCP} we put all
pieces together.

\subsection{Orientation preserving case of the TCP: periodicity boxes and building
conjugators\label{ssec:periodicity-boxes}}

We now deal with the equation $z=g^{-1}yg$ for $y,z\in \EP_2$ and $g\in F$. The argument
will make use of techniques and statements in~\cite{matucci5} and refer often to that
paper.

Subsection~4.1 in~\cite{matucci5} shows that, if $z=g^{-1}yg$ with $y,z,g \in \PL_2(I)$,
then there exists $\varepsilon
>0$ depending only on $y$ and $z$ such that $g$ is linear inside $[0,\varepsilon]^2$; the box
$[0,\varepsilon]^2$ is called an \emph{initial linearity box}. The goal of this section is
to show an analog of this result inside suitable boxes $(-\infty,L]^2$ and $[R,\infty)^2$
where $y,z \in \EP_2$ are periodic.

The following is a first necessary condition for two maps to be conjugate to each other.

\lemma{\label{thm:identical-at-infinity}
Let $y,z \in \EP_2$ be such that $y\sim_F z$. Then there exist two numbers $L,R \in
\mathbb{R}$ such that $y(t)=z(t)$ for all $t \in (-\infty,L] \cup [R,\infty)$.
}

\begin{proof} Let $g\in F$ be such that $g^{-1}y g=z$. For $t$ negative sufficiently large, we have $g(t)=t+m_-$, and so
 $$
z(t)=g^{-1} y g(t)= g^{-1}y(t+m_-)=g^{-1}(y(t)+m_-)=y(t)+m_- -m_-=y(t).
 $$
Similarly for $t$ positive sufficiently large.
\end{proof}

We move on to prove the existence of periodicity boxes.

\lemmaname{Initial and final periodicity boxes}{\label{ifpb}
For every pair of elements $y,z \in \EP_2^>(-\infty, p)$ (with $-\infty <p\leqslant
+\infty$), there exists a computable constant $L\in \mathbb{R}$ (depending only on $y$ and
$z$) such that every conjugator $g\in F$ between $y$ and $z$ must act as
a translation inside the \emph{initial periodicity box} $(-\infty,L]^2$. Similarly, for
every pair of elements $y,z \in \EP_2^>(p, +\infty)$ (with $-\infty \leqslant p<+\infty$)
and a \emph{final periodicity box} $[R, +\infty)^2$.

The exact same statement is true replacing $\EP_2^>$ to $\EP_2^<$.
}

\begin{proof}
If $y$ and $z$ are not equal for $t$ positive and negative sufficiently large then, by
Lemma~\ref{thm:identical-at-infinity}, there is no possible conjugator $g\in F$ and there
is nothing to prove. So assume they are and consider a negative sufficiently large $L\in
\mathbb{R}$ such that $y(t)=z(t)$ and $y(t-1)=y(t)-1$ (and so, $z(t-1)=z(t)-1$), for every
$t\leqslant L$ (clearly, such an $L$ is computable). We claim that every possible $g\in F$
satisfying $g^{-1}yg=z$ must be a translation for $t\leqslant L$. By the symmetry of $y$
and $z$ in the definition of $L$ and up to writing the conjugacy relation as
$(g^{-1})^{-1}zg^{-1}=y$ (which changes the conjugator from $g$ to $g^{-1}$), we can assume
that $g$ has non-positive translation at $-\infty$ (i.e. $g(t)=t+m_-$ for negative
sufficiently large $t$, and with $m_-\leqslant 0$).

Assume, by contradiction, that $g$ is not a translation map in $(-\infty,L]$. Then, there
is $\lambda<L$ such that
 $$
g(t)=\begin{cases} t+m_- & t \leqslant \lambda \\ \alpha(t-\lambda) + \lambda +m_- & \lambda
\leqslant t<\mu \end{cases}
 $$
for some suitable real numbers $\alpha \neq 1$, $\lambda <\mu <L$. Since $z$ is increasing
and strictly above the diagonal $\id(t)=t$, we can choose $r<\lambda <L$ such that $\lambda
<z(r)<\mu <L$. By our choice of $r$, we have $y(r)=z(r)$, $y(t-1)=y(t)-1$ and
$z(t-1)=z(t)-1$ for all $t\leqslant r$. Moreover, since $gz(t)=yg(t)$ for all $t \in
\mathbb{R}$, we have
 $$
\alpha(z(r)-\lambda) + \lambda +m_- =gz(r) =yg(r) =y(r+m_-) =y(r)+m_- =z(r)+m_-.
 $$
Rearranging the terms, we have
 $$
\alpha(z(r)-\lambda)=z(r)-\lambda
 $$
and, since $z(r)-\lambda >0$, we get $\alpha=1$, a contradiction. Hence, $g(t)=t+m_-$ for
every $t\leqslant L$ as claimed.

The symmetric argument gives a constant $R$ establishing the final periodicity box
$[R,+\infty)^2$.

If $y,z\in \EP_2^<$, then we apply the previous argument to $y^{-1},z^{-1}$ and derive the
same conclusion.
\end{proof}

\remark{Note that, in the previous lemma, the constants $L$ and $R$ depend on $y$ and $z$
but not on the conjugator $g$. This will be crucial later.}

We observe that the results of Subsection 4.2 in \cite{matucci5} and their proofs follow
word-by-word in our generalized setting, and hence we do not reprove them. We restate Lemma
4.6 in \cite{matucci5} to give an example of how results appear in this context.

%
%
%
%
%

\lemma{ Let $z\in \EP_2^<$. Let
$C_{F}(z)=C_{\PL_2(\mathbb{R})}(z)\cap F$ be the set of elements in $F$ commuting with $z$.
Then the map $\varphi_z:C_{F}(z) \to \mathbb{Z}$ defined by
 $$
\varphi_z(g) = \lim_{t \to -\infty}g(t)-t
 $$
is an injective group homomorphism. A similar statement is true for $\EP_2^>$.}



Subsection 4.2 in \cite{matucci5} shows how to build a candidate conjugator $g$ between any
two elements of $F$ after we have chosen the initial slope of $g$.

A \emph{unique candidate conjugator} $g$ between $y$ and $z$ with a given initial slope
$q$, if it exists, is the unique function that one needs to test as a conjugator of $y$ and
$z$ with initial slope $q$: if $g$ fails to satisfy $g^{-1}yg=z$, then there is no
conjugator of $y$ and $z$ with initial slope $q$. The proof of Corollary 4.12 in
\cite{matucci5} can be lifted verbatim and so we only restate it in our new case.


\theoremname{Explicit conjugator}{ \label{thm:explicit-conjugator} Let $y,z \in
\EP_2^<$. Suppose there exist $L<R$ such that $y$ and $z$ coincide and are periodic on $(-\infty,L] \cup [R,+\infty)$,
so that $(-\infty,L]^2$ is the initial periodicity box.
Let $\ell \in \mathbb{Z}_{<0}$.
\begin{enumerate}
\item Let $g_0 \in F$ be a map which is affine inside $(-\infty,L)^2$ and such that
    $\lim_{t \to -\infty}g_0(t)-t=q$. Then the unique conjugator $\wh{g} \in \PL_2(\R)$
    between $y$ and $z$, which is affine inside $(-\infty,L)^2$ and such that $\lim_{t
    \to -\infty}\wh{g}(t)-t=\ell$ is defined pointwise by
$$
\wh{g}(t)=\lim_{r\to +\infty} y^{-r}g_0z^r(t).
$$
Moreover, the map $\wh{g}$ is recursively constructible and $y$ and $z$ are always conjugate
in $\PL_2(\mathbb{R})$ via $\wh{g}$.
\item There exists an algorithm to decide whether or not there is $g \in F$ such that
    $\lim_{t \to -\infty}g(t)-t=\ell$ and $g^{-1}yg=z$.
\end{enumerate}
}

The above result has been stated, for simplicity, for two functions $y,z \in \EP_2^<$.
However, the same result can be stated for $y,z \in \PL_2^<([p_1,p_2])$ for any $p_1, p_2
\in \mathbb{Q}$, or for $y,z \in \EP_2^<(p,+\infty)$.

\remark{The results of this subsection do not involve dyadic rationals and slopes that are
powers of 2 and are, in fact, true for other classes of groups without restrictions on the
breakpoints and the slopes (for example $\PL_+(\mathbb{R})$, the
Bieri-Thompson-Stein-Strebel groups in $\mathbb{R}$ and the corresponding subgroups with
eventually periodic tails). See \cite{matucci5} for more details.}

\subsection{Orientation preserving case of the TCP: fixed points\label{ssec:fixed-points}}

The goal of this Subsection is to reduce to the case where the sets $\partial \Fix (y)$ and
$\partial \Fix (z)$ do coincide. Up to suitable special cases, this will allow us to reduce
to looking for potential conjugators $g\in F$ such that $\partial \Fix (y)=\partial
\Fix(z)\subseteq \Fix (g)$, thus restricting ourselves to studying conjugacy among the
corresponding intervals of $y$ and $z$ between any two consecutive points $p$ and $q$ of
$\partial \Fix(y)=\partial \Fix(z)$. On each such interval $y$ (and $z$) is either the
identity, or has no fixed points apart from $p$ and $q$ and so they belong to either
$\EP_2^<(p,q)$ or $\EP_2^>(p,q)$.

Note that the sets $\partial \Fix(y)$ and $\partial \Fix(z)$ are discrete subsets of
$\mathbb{Q}$, and their intersections with any finite interval $[L, R]$ are easily
computable by just solving finitely many systems of linear equations. An apparent technical
difficulty is that, since $y,z \in \EP_2$, the full sets $\partial \Fix(y)$ and $\partial
\Fix(z)$ may be infinite; however, due to the periodicity, they are controlled by finite
sets.

\proposition{\label{thm:identify-fixed-points}
There is an algorithm which, given $y,z\in \EP_2$ being equal for $t$ negative sufficiently
large and for $t$ positive sufficiently large, decides whether or not there exists some
$g\in F$ such that $\partial \Fix(y)=g(\partial \Fix(z))$ and, in the affirmative case, it
constructs such a $g$.
}

\begin{proof}
For the given $y,z$ we can easily compute constants $L<R$ such that, for all $t\in
(-\infty, L]$, $y(t)=z(t)$ and $y(t-1)=y(t)-1$, and such that, for all $t\in [R, +\infty)$,
$y(t)=z(t)$ and $y(t+1)=y(t)+1$. Moving $L$ down and/or $R$ up if necessary, we can also
assume that if $\partial \Fix(y) \neq \emptyset$ then it has at least one point in $[L, R)$
(and similarly for $z$).

Now compute the finite sets of rational numbers $\partial \Fix(z)\cap [L, R)$, $\partial
\Fix(y)\cap [L, R)$, $\partial \Fix(z)\cap [L-1, L)=\partial \Fix(y)\cap [L-1, L)$,  and
$\partial \Fix(z)\cap [R, R+1)=\partial \Fix(y)\cap [R, R+1)$; let $p,q,m,n\geqslant 0$ be
their cardinals, respectively. By the periodicity of $y$ and $z$ outside $[L, R]$, these
constitute full information about $\partial \Fix(y)$ and $\partial \Fix(z)$. Up to
switching $y$ with $z$, we may assume that $p\leqslant q$.

Clearly, $m=0$ if and only if $\partial \Fix(y)$ and $\partial \Fix(z)$ have a minimum
element (as opposed to having infinitely many points approaching $-\infty$). Similarly,
$n=0$ if and only if $\partial \Fix(y)$ and $\partial \Fix(z)$ have a maximum element.

If either $\partial \Fix(y)$ or $\partial \Fix(z)$ is empty then there is nothing to prove.
Assume $\partial \Fix(y)\neq \emptyset \neq \partial \Fix(z)$, i.e. $1\leqslant p\leqslant
q$. We denote by $a_0$ (respectively, $b_0$) the smallest element in $\partial \Fix(z)\cap
[L,R)$ (respectively $\partial \Fix(y)\cap [L,R)$) and we use it to enumerate in an order
preserving way all the elements of the discrete set $\partial \Fix(z)$ (respectively,
$\partial \Fix(y)$) as $a_i$ (respectively, $b_i$); the index $i$ will run over a finite,
infinite or bi-infinite subset of $\mathbb{Z}$ depending on whether or not $m$ (and/or $n$)
is zero. With this definition, $\partial \Fix(z)\cap [L, R)=\{ a_0<a_1<\cdots <a_{p-1} \}$
and $\partial \Fix(y)\cap [L, R)=\{ b_0<b_1<\cdots <b_{q-1} \}$.

Note that any $g\in F$ satisfying $\partial \Fix(y)=g(\partial \Fix(z))$ must map all the
$a_i$'s bijectively to all the $b_i$'s. In particular, if $m=0$ then $a_0$ must be mapped
to $b_0$, and if $n=0$ then $a_{p-1}$ must be mapped to $b_{q-1}$ (and so $a_0$ to
$b_{q-p}$). Hence, in the special case that either $m=0$ or $n=0$, the following claim
completes the proof.

\bigskip
\noindent \emph{Claim 1:} For every $b_i\in \partial \Fix(y)$, we can algorithmically
decide whether or not there exists some $g\in F$ such that $\partial \Fix(y)=g(\partial
\Fix(z))$ and $g(a_0)=b_i$ and, in the affirmative case, the algorithm constructs one
explicitly.

\bigskip

The remaining case to study is when $m\neq 0\neq n$, so that $a_0$ potentially could be
sent to any of the $b_i$'s by the map $g$. Let $\ell=\mathrm{lcm}(m,n)$ and let $[L-\ell/m,
L)$ be the smallest interval to the left of $L$ to contain $\ell$ points of $\partial
\Fix(z)$. Similarly, let $[R, R+\ell/n)$ be the corresponding interval to the right of $R$.
Consider the following two finite sets:
 $$
\begin{array}{c}
A:=\partial \Fix(z) \cap \left[L-\frac{2\ell}{m}, R+\frac{2\ell}{n} \right), \\ \\
B:=\partial \Fix(y) \cap \left[L-\frac{2\ell}{m}, R+\frac{2\ell}{n} \right),
\end{array}
 $$
and let $s_0$ be the rightmost point of $\partial \Fix(z) \cap \left[L-\frac{2\ell}{m},
L-\frac{\ell}{m} \right)$, and let $t_0$ be the leftmost point of $\partial \Fix(z) \cap
\left[R+\frac{\ell}{n}, R+\frac{2\ell}{n} \right)$. We compute $A$, $B$, $s_0$ and $t_0$
explicitly.

\bigskip
\noindent \emph{Claim 2:} Suppose there exists a map $g\in F$ such that $\partial
\Fix(y)=g(\partial \Fix(z))$ and $g(s_0)\in \left[R+\frac{k\ell}{n},
R+\frac{(k+1)\ell}{n}\right)$ for $k\geqslant 2$; then, there exists a $g'\in F$ such that
$\partial \Fix(y)=g'(\partial \Fix(z))$ and $g'(s_0)\in \left[R+\frac{(k-1)\ell}{n},
R+\frac{k\ell}{n}\right)$. Similarly, if there exists $g\in F$ such that $\partial
\Fix(y)=g(\partial \Fix(z))$ and $g(t_0)\in \left[L-\frac{(k+1)\ell}{n},
L-\frac{k\ell}{n}\right)$ for some $k\geqslant 2$, then there exists a $g'\in F$ such that
$\partial \Fix(y)=g'(\partial \Fix(z))$ and $g'(t_0)\in \left[L-\frac{k\ell}{n},
L-\frac{(k-1)\ell}{n}\right)$.

\bigskip
With the help of Claim~2 we can complete the proof in the following way. Suppose there
exists $g\in F$ such that $\partial \Fix(y)=g(\partial \Fix(z))$. Since $p\leqslant q$ it
cannot simultaneously happen that $g(s_0)<s_0$ and $t_0<g(t_0)$. Hence either $s_0\leqslant
g(s_0)$ or $g(t_0)\leqslant t_0$ and, in either case, a repeated application of Claim~2
implies the existence of $g'\in F$ such that $\partial \Fix(y)=g'(\partial \Fix(z))$ and
$g'(A)\cap B\neq \emptyset$. This gives finitely many possibilities for $g'(a_0)$ and so,
applying Claim~1 finitely many times we can decide whether or not there exists a $g\in F$
satisfying $\partial \Fix(y)=g(\partial \Fix(z))$.

Hence, it only remains to prove the above two claims.

\bigskip
\begin{proof}[Proof of Claim~1]
We will distinguish four cases.

\bigskip
\noindent \textbf{Case 1: $m=0$ and $n=0$.} In this case, $\partial \Fix(z)=\{
a_0<a_1<\cdots <a_{p-1} \}$ and $\partial \Fix(y)=\{ b_0<b_1<\cdots <b_{q-1} \}$ and,
clearly, $p=q$ and $g(a_0)=b_0$ are necessary conditions for such a $g$ to exist. If both
conditions hold, then Proposition~\ref{thm:rationals-coincide} makes the decision for us.

\bigskip
\noindent \textbf{Case 2: $m\geqslant 1$ and $n=0$.} This case is entirely symmetric to the
next one.

\bigskip
\noindent \textbf{Case 3: $m=0$ and $n\geqslant 1$.} In this case, $\partial \Fix(z)$ and
$\partial \Fix(y)$ both have first elements $a_0$ and $b_0$ and infinitely many points
approaching $+\infty$. As in case~1, $g(a_0)=b_0$ is a necessary condition for such a $g$
to exist.

We have $\partial \Fix(z)\cap [R, R+1)=\{ a_{p}<a_{p+1}<\cdots <a_{p+(n-1)} \}$ and that
the elements in $\partial \Fix(z) \cap [R+1,+\infty)$ are integral translations of these:
for every $j\geqslant 0$, write $j=\lambda n+\mu$ with $\lambda,\mu\geqslant 0$ integers
and $\mu =0,\ldots, n-1$, and we have $a_{p+j}=\lambda+a_{p+\mu}$. A similar argument for
$y$ yields that $\partial \Fix(y)\cap [R, R+1)=\{ b_{q}<b_{q+1}<\cdots <b_{q+(n-1)} \}$ and
that, for every $j\geqslant q$, we have $b_{q+j}=\lambda+b_{q+\mu}$. Moreover, from
$a_p=b_q$ on, the two sequences coincide, i.e., for every $j\geqslant 0$,
 $$
\lambda+a_{p+\mu}=a_{p+j}=b_{q+j}=\lambda+b_{q+\mu}.
 $$

Now if some $g\in F$ satisfies $g(\partial \Fix(z))=\partial \Fix(y)$, it must apply the
points in an order preserving way, starting from the smallest ones, that is, $g(a_k)=b_k$
for any integer $k$. In particular, for $k\geqslant q\geqslant p$, we have
 $$
g(\lambda_1+a_{p+\mu_1})=g(a_{p+(k-p)})=g(a_k)=b_k=b_{q+(k-q)}=\lambda_2+b_{q+\mu_2},
 $$
where $k-p=\lambda_1 n+\mu_1$ and $k-q=\lambda_2 n+\mu_2$. Since $g$ is of the form
$g(t)=t+m_+$ with $m_+\in \mathbb{Z}$ for $t$ positive sufficiently large then, for large
enough $k$, the above equation tells us that
 $$
\lambda_1+a_{p+\mu_1}+m_+=g(\lambda_1+a_{p+\mu_1})=\lambda_2+b_{q+\mu_2}.
 $$
Therefore, $a_{p+\mu_1}-b_{q+\mu_2}=b_{q+\mu_1}-b_{q+\mu_2}$ must be an integer and so,
$\mu_1=\mu_2$, which means that $k-p$ and $k-q$ are congruent modulo $n$, i.e. $q-p$ is
multiple of $n$.

Assume then this necessary condition, $q-p=\lambda n$ with $\lambda \in \mathbb{Z}$, and
apply Proposition~\ref{thm:rationals-coincide}~(2) to the sequences $a_0<\cdots
<a_{p+\lambda n-1}$ and $b_0<\cdots <b_{q-1}$ (both with $q$ points). If there is no $g\in
F$ sending the first list to the second then there is no $g$ such that $\partial \Fix(y)
=g(\partial \Fix(z))$ and we are done. Otherwise, we get a $g$ matching these first $q$
points, $g(a_0)=b_0, \ldots ,g(a_{p+\lambda n-1})=b_{q-1}$, and, after a final small
modification, we will see that it automatically matches the rest.

Choose two dyadic numbers $a_{p+\lambda n-1}<\alpha <\beta <a_{p+\lambda n}$, choose $h\in
F$ such that $h(\alpha )=g(\alpha)$ and $h(\beta )=\beta -\lambda$ (such an $h$ exists and
is effectively computable by Proposition~\ref{thm:rationals-coincide}~(2)), and let us
consider the following map:
 $$
\widetilde{g}(t)=\begin{cases} g(t) & t \leqslant \alpha \\ h(t) & \alpha \leqslant t\leqslant \beta
\\ t-\lambda & \beta \leqslant t.\end{cases}
 $$
By construction, $\widetilde{g}$ is continuous, piecewise linear with dyadic breakpoints,
and all slopes are powers of 2; furthermore $g\in F$ and $\widetilde{g}$ is an integral
translation for $t\geqslant \beta$ so, $\widetilde{g}\in F$. On the other hand,
 $$
\partial \Fix(y)\cap [L, b_{q-1}]=g(\partial \Fix(z)\cap [L,a_{p+\lambda
n-1}])=\widetilde{g}(\partial \Fix(z)\cap [L,a_{p+\lambda n-1}]),
 $$
and
 $$
\partial \Fix(y)\cap [b_{q}, +\infty)=\{ b_q, b_{q+1},\ldots \}=\{ a_{p+\lambda n}-\lambda,
a_{p+\lambda n+1}-\lambda, \ldots \} =
 $$
 $$
=\widetilde{g}(\{ a_{p+\lambda n}, a_{p+\lambda n+1}, \ldots \})=\widetilde{g}(\partial \Fix(z)\cap
[a_{p+\lambda n}, +\infty )).
 $$
Hence, $\partial \Fix(y)=\widetilde{g}(\partial \Fix(z))$ and we are done.

\bigskip
\noindent \textbf{Case 4: $m\geqslant 1$ and $n\geqslant 1$.} The argument in this case is
similar to that of case~3 but repeated twice, up and down (and with no restriction for
$b_i$ because we have both infinitely many fixed points bigger and smaller than $b_i$).

Following the notation above, the $m$ fixed points from $\partial \Fix(z)\cap [L-1,
L)=\partial \Fix(y)\cap [L-1, L)$ are labeled and ordered as $a_{-m}<\cdots <a_{-1}$ and
$b_{-m}<\cdots <b_{-1}$ (hence, $a_{-j}=b_{-j}$ for $j=1,\ldots ,m$). The elements from
$\partial \Fix(z) \cap (-\infty,L-1)$ and $\partial \Fix(y) \cap (-\infty,L-1)$ are their
integral translations to the left.

Now if some $g\in F$ satisfies $g(\partial \Fix(z))=\partial \Fix(y)$ and $g(a_0)=b_i$, it
must send the points $a_j$ to the $b_j$ in an order preserving way starting from
$g(a_0)=b_i$, both up and down. Hence, two arguments exactly like in the previous case give
us two necessary congruences among $p,q$ and $i$, modulo $n$ (close to $+\infty$) and
modulo $m$ (close to $-\infty$). If one of them fails, then there is no such $g$ and we are
done. If both are satisfied, then apply Proposition~\ref{thm:rationals-coincide}~(2) to a
long enough tuple of $a_j$'s and $b_j$'s: a negative answer tells us there is no such $g\in
F$, and a positive answer provides a $g\in F$ which, after two local modifications like in
the previous case (one close to $+\infty$ and the other close to $-\infty$), will finally
give us a $g'\in F$ such that $g'(\partial \Fix(z))=\partial \Fix(y)$, and $g'(a_0)=b_i$.

This completes the proof of Claim~1.
\end{proof}

\begin{proof}[Proof of Claim 2]
We will prove the first part of the claim; the symmetric argument for the second part is
left to the reader.

Assume the existence of $g\in F$ such that $g(\partial \Fix(z))=\partial \Fix(y)$ and
$g(s_0 )\in \left[R+\frac{k\ell}{n}, R+\frac{(k+1)\ell}{n}\right)$ for $k\geqslant 2$. To
push $g(s_0)$ down, let us define the \emph{reduction} map $g_-$ by
$$
g_-(t)=
\begin{cases}
g(t-\frac{\ell}{m}) & t<s_0 \\
g(t)-\frac{\ell}{n} & t \ge s_0.
\end{cases}
$$
To understand the map $g_-$, note that its graphical representation can be obtained from
that of $g$ by performing the following operation:
remove the graph within $[s_0-\ell/m,s_0]$, translate the graph of $g$ defined on
$[s_0,+\infty)$ by the vector $(0,-\ell/m)$ and translate the graph of $g$ defined on
$(-\infty,s_0-\ell/m]$ by the vector $(\ell/m,0)$. Hence, $g_-$ is the same as $g$ avoiding
the piece over the interval $[s_0-\ell/m, s_0]$.

It is obvious that the two parts of $g_-$ to the left and to the right of $s_0$ are both
continuous, increasing, piecewise linear, with dyadic breakpoints, with slopes being powers
of two, and being eventually translations (near $-\infty$ and $+\infty$, respectively). To
check whether $g_-$ is in $F$ it only remains to analyze what happens around the point
$s_0$.

First of all, $g_-$ is continuous at $s_0$: observe that $s_0-\frac{\ell}{m}\in
\partial \Fix(z)$ is exactly $\ell$ points to the left of $s_0$ in the discrete set
$\partial \Fix(z)$; since $g(\partial \Fix(z))=\partial \Fix(y)$ and $g$ is an increasing
function, $g(s_0-\frac{\ell}{m})$ must be exactly $\ell$ points to the left of $g(s_0)$ in
the discrete set $\partial \Fix(y)$ that is, $g(s_0-\frac{\ell}{m})=g(s_0)-\frac{\ell}{n}$.

Unfortunately, the slopes of $g_-$ to the left and to the right of $s_0$ (i.e. the slopes
of $g$ to the left of $s_0-\ell/m$ and to the right of $s_0$) may be different; and $s_0$
may not be a dyadic rational number. If these two facts happen simultaneously then $g_-$
will not an element of $F$ because of having a breakpoint at a non-dyadic point, namely
$s_0$. This technical difficulty will be fixed later by modifying the map $g_-$ in a
suitably small neighborhood of $s_0$.

Before doing this, let us check that $g_-$ fulfils our requirement. Since
$g(s_0-\frac{\ell}{m})=g(s_0)-\frac{\ell}{n}$, the hypothesis $g(\partial \Fix(z))=\partial
\Fix(y)$ implies that
 $$
g_-(\partial \Fix(z)\cap (-\infty, s_0])=g\left(\partial \Fix(z)\cap \left(-\infty,
s_0-\frac{\ell}{m}\right]\right)=
 $$
 $$
=\partial \Fix(y)\cap \left(-\infty, g(s_0)-\frac{\ell}{n}\right],
 $$
and $g_-(s_0)=g(s_0-\frac{\ell}{m})=g(s_0)- \frac{\ell}{n}$, and
 $$
g_-(\partial \Fix(z)\cap [s_0, +\infty))=g(\partial \Fix(z)\cap [s_0, +\infty ))-
\frac{\ell}{n} =\partial \Fix(y)\cap \left[g(s_0)-\frac{\ell}{n}, +\infty \right).
 $$
Hence, $g_-(\partial \Fix(z))=\partial \Fix(y)$ and $g_-(s_0)=g(s_0)-\frac{\ell}{n}\in
\left[R+\frac{(k-1)\ell}{n}, R+\frac{k\ell}{n}\right)$, as we wanted.

To complete the proof of Claim~1 we must be able to fix the above technical problem, by
modifying $g_-$ in such a way that the resulting map belongs to $F$, but not changing the
image of any point in $\partial \Fix(z)$; this will be achieved by changing $g_-$ only in a
small enough neighborhood of $s_0$ not containing any other point of $\partial \Fix(z)$
(and, of course, not changing the image of $s_0$ itself).

Let $\alpha_1$ be a dyadic point found strictly between $\alpha_2:=s_0$ and the point of
$A$ immediately to the left of $s_0$; and let $\alpha_3$ be a dyadic point found strictly
between $\alpha_2:=s_0$ and the point of $A$ immediately to the right of $s_0$. Now
consider the points
 $$
\beta_1:=g_-(\alpha_1)=g\left(\alpha_1-\frac{\ell}{m} \right),
 $$
 $$
\beta_2:=g_-(\alpha_2)=g\left(\alpha_2-\frac{\ell}{m} \right)=g(\alpha_2)-\frac{\ell}{n},
 $$
 $$
\beta_3:=g_-(\alpha_3)=g(\alpha_3)-\frac{\ell}{n}.
 $$
Since $\alpha_1<\alpha_2<\alpha_3$ and $\beta_1<\beta_2<\beta_3$ are rational points such
that, for every $i=1,2,3$, $\beta_i$ is the image of $\alpha_i$ by some element in $F$,
then we can apply Proposition~\ref{thm:rationals-coincide}~(2) and construct a function
$h\in F$ such that $\beta_i=h(\alpha_i)$. Finally, define
 $$
g'(t)=\begin{cases} h(t) & t\in [\alpha_1,\alpha_3] \\ g_-(t) & t\not \in [\alpha_1,\alpha_3].
\end{cases}
 $$
Clearly, $g'\in F$, $g'(s_0)=g_-(s_0)$ and $g'(\partial \Fix(z))=g_-(\partial
\Fix(z))=\partial \Fix(y)$. This completes the proof of Claim~2.
\end{proof}
This finishes the proof of Proposition~\ref{thm:identify-fixed-points}.
\end{proof}

\lemma{\label{thm:RTCP} The decidability of the following two problems is equivalent:
\begin{enumerate}
\item[(TCP)] For any two $y,z\in \EP_2$ we can determine whether or not there is $g\in F$
    such that $g^{-1}yg=z$.
\item[(RTCP)] For any two $y,z\in \EP_2$ such that $\partial \Fix(y)=\partial \Fix(z)$ we
    can determine, whether or not there is $g\in F$ such that $g^{-1}yg=z$.
\end{enumerate}
}

\begin{proof}
Obviously, if (TCP) is decidable, then (RTCP) is decidable. Assume now that (RTCP) is
decidable. By the discussion at the beginning of this subsection, if $y$ and $z$ are
conjugate via $g \in F$, then $\partial \Fix(y)=g(\partial \Fix(z))$. By
Theorem~\ref{thm:identify-fixed-points} we can decide whether or not there is a map $g \in
F$ such that $\partial \Fix(y)=g(\partial \Fix(z))$. If there is no such map, then $y$ and
$z$ are not conjugate. If there is such a $g\in F$ (and in this case
Theorem~\ref{thm:identify-fixed-points} constructs it) then $\partial
(\Fix(gzg^{-1}))=g(\partial \Fix (z))=\partial \Fix (y)$ and we can apply (RTCP) to the two
maps $y$ and $gzg^{-1}$ to detect whether or not they are conjugate. Obviously, this is the
same decision as the one we are interested in.
\end{proof}

By Lemma \ref{thm:RTCP} we can restrict our focus to studying (RTCP).

\subsection{Orientation preserving case of the TCP: Reducing the problem to squares.
\label{ssec:reducing-to-squares}}

We can make $\partial \Fix(y)=\partial \Fix(z)$ as done in
Proposition~\ref{thm:identify-fixed-points}. If $\partial \Fix(y)=\partial \Fix(z) =
\emptyset$ we defer the discussion to Subsection \ref{sec:rescaling-the-circle}. On the
other hand, if $\partial \Fix(y)=\partial \Fix(z) \neq \emptyset$ and $g\in F$ is a
conjugator between $y$ and $z$, the only thing we can say is that $g$ acts on $\partial
\Fix(y)$ in an order preserving way. There are two possibilities:
\begin{enumerate}
\item $\Fix(g) \neq \emptyset$.
\item $\Fix(g) = \emptyset$. We can assume that $g \in F^>$.
\end{enumerate}
Case (2) can indeed happen as is shown by the following example: take $y=z$ to be a
non-trivial periodic function of period $1$ with fixed points. Then the map $g(t)=t+1 \in
F^>$ is a conjugator for $y$ and $z$ having no fixed points.


We need to find if there is a conjugator $g$ between $y$ and $z$ such that $g \in F^>$. 
We can assume $y \ne \id \ne z$, otherwise our analysis becomes trivial.
We can write the supports $\supp(y)=\supp(z)$ as the union of the family $\{I _j\}$ of
(possibly unbounded) intervals on which $y$ and $z$ have no fixed points ordered so that
$I_j$ is to the left of $I_{j+1}$, for every $j$. 
If this family were finite, since we are assuming $\partial \Fix(y)=\partial \Fix(z) \neq \emptyset$,
then it means that $\partial \Fix(y)$ is finite and so $g$ must fix the smallest element in $\partial \Fix(y)$
since $g$ is order preserving, hence $\Fix(g) \ne \emptyset$ and this would not be
the case that we are studying now.  Thus we must study the case of
the following proposition.

\noindent \proposition{Let $y,z \in \EP_2$ be such that $\Fix(y)=\Fix(z)$ and that
$\partial \Fix(y)$ has infinitely many points. Then there are only finitely many candidate
conjugators $g \in F^>$.}

\begin{proof}
We give out only some relevant details of how to prove this proposition. This entails
generalizations of many results of this paper and of Kassabov-Matucci ~\cite{matucci5} and
so we only explain how to achieve them. The main point here is noticing that we can develop
a Stair Algorithm and bounding initial slopes of $g \in F$, even if at $-\infty$ the
functions $y,z$ have no initial slope.

By hypothesis, $\{I_j\}$ has infinitely many intervals and so $g$ ``shifts'' them, that is
$g(I_j)=I_{j+k}$, for some fixed $k$. Let $t_j$ be the left endpoint of $I_j$. We make a
series of observations:
\begin{enumerate}
\item We can build candidate conjugators (Theorem \ref{thm:explicit-conjugator}) on each
    $I_j$, given a fixed initial slope at $t_j$,
\item The initial slope of $z$ on $I_j$ coincides with the initial slope of $y$ in the
    image interval $g(I_j)$,
\item There is an ``initial'' box for $g$ in $I_j$,
\item We can bound the ``initial'' slopes of $g$ on $I_j$,
\item We can bound the initial slope of $g$ at $-\infty$.
\end{enumerate}

(1) and (2) are a straightforward calculation. (3) is a verbatim rewriting of the proof of
Lemma 4.2 in ~\cite{matucci5}.

(4) A standard trick from ~\cite{matucci5} is observing that
\[
z=g^{-1}yg= g^{-1}y^{-r} y y^rg
\]
and so the slope of $y^r g$ at $t_j$ is $(y'(t_{i+k}))^r g'(t_i)$ and $y^r g$ is a
conjugator for $y$ and $z$ on $I_j$. On each $I_j$ there are only finitely many slopes for
$g'(t_i^+)$ to be tested and on each one, we apply Theorem \ref{thm:explicit-conjugator} to
build candidate conjugators that we can test.

(5) Recall that a candidate conjugator $g$ pushes all the intervals in $\supp(y)$ in the
same direction by the ``same amount of intervals in $\supp(y)$''. In particular, the
initial slope of $g$ determines the number $k$ such that $g(I_j)=I_{j+k}$ for every $j$.

We use ideas similar to Claim 2 in Proposition \ref{thm:identify-fixed-points}. Let us call
$J_L$ the left open interval on which $y=z$ and they are periodic. A similar definition can
be made for $J_R$. Let $J_C=\mathbb{R} \setminus (J_L \cup J_R)$ the remaining central
piece. Assume that there is a conjugator $g$ between $y$ and $z$ which sends and interval
$I_j$ inside $J_L$ to an interval $I_{j+k+1}$ with the requirement that $I_{j+k}$ is
entirely contained into $J_R$. Using ideas similar to Claim 2 in Proposition
\ref{thm:identify-fixed-points} one can create a new conjugator $\overline{g}$ such that
$\overline{g}(I_j)=I_{j+k}$.

Therefore, similarly to Claim 2 in Proposition \ref{thm:identify-fixed-points}, this allows
us to reduce the study to only finitely many candidate conjugators where $g(J_C) \cap J_C
\ne \emptyset$ or where the rightmost interval $I_j$ inside $J_L$ goes to the leftmost
interval $I_s$ of $J_R$ (or viceversa). This argument reduces the number of initial slopes
of $g$ to be tested.

To conclude we observe that there are only finitely many slopes for $g$ at $-\infty$ and
finitely many ``initial'' slopes for $g$ on the finitely many intervals $I_j$ contained in
$J_C$ and then we can apply Theorem \ref{thm:explicit-conjugator} on each of these
intervals building finitely many candidate conjugators $g \in F^>$ which we can then test
one by one.
\end{proof}

The previous result allows one to restrict to the case of looking for conjugators $g$ with
fixed points.

\lemma{\label{tim:fix-of-z-is-inside-fix-g} Let $y,z \in \EP_2$ be such that
$\Fix(y)=\Fix(z) \ne \emptyset$ and let $g$ be a conjugator between $y$ and $z$ such that
$\Fix(g) \ne \emptyset$. Then $\Fix(z) \subseteq \Fix(g)$.}

\begin{proof}
Let $a \in \Fix(g)$ and let $b$ be the the smallest point of $\partial \Fix(z)$ such that
$a<b$. Since $g$ fixes $\Fix(z)$ set wise and is order-preserving, then $g(b)$ must also be
the smallest point of $\partial \Fix(z)$ such $g(b)>a$, therefore $g(b)=b$ and so $g$ must
fix all of $\Fix(z)$ pointwise.
\end{proof}

We need to show that (RTCP) of Lemma \ref{thm:RTCP} is decidable. Lemma
\ref{tim:fix-of-z-is-inside-fix-g} tells us that we can restrict ourselves to solve the
problem inside the closed intervals of $\Fix(y)=\Fix(z)$.

As was done in \cite{matucci5} we observe that if $p \in \partial \Fix(y)$
is a non-dyadic rational point and $g$ is a conjugator between $y$ and $z$, then
$g'(p^-)=g'(p^+)$ or, in other words, the slope of $g$ at one side of $p$ is completely determined
by the slope on its other side. This implies that the important points of $\partial \Fix(y)$
are the dyadic rational ones (if they exist) as they are the ones where $g$ has freedom
to have different slopes on the two sides and therefore the conjugator that we are attempting to build
can be constructed by by gluing two conjugators on the two sides of a dyadic rational point of $\partial \Fix(y)$.
In the case that $\partial \Fix(y)$ had no dyadic rational points, then we can compute a conjugator at a point
$p \in \partial \Fix(y)$ and this uniquely determines the conjugator on the entire real line.
Otherwise, there are dyadic rational points in $\partial \Fix(y)$ and we argue as following.

Let $L<R$ are two integers chosen
so that $y$ and $z$ coincide and are periodic inside $(-\infty,L] \cup [R,+\infty)$. The case
when $\partial \Fix(y) \cap [L,R]$ contains no dyadic rational point is dealt with as above. Similarly, if
there is only one dyadic point inside $\partial \Fix(y) \cap [L,R]$, then we have two instances of the previous case on
the two sides of the dyadic point. Otherwise, we choose $p_1,p_2$ with the property of being dyadic and consecutive inside $\partial \Fix(y)$
and such that $[p_1,p_2] \subseteq [L,R]$. With these provisions, we can use the
solution of the standard conjugacy problem inside $\PL_2([p_1,p_2])$ using the techniques
in \cite{matucci5}. 
If there is no conjugator on any of those intervals, then $y$ and $z$
cannot be conjugate. Otherwise, we can glue the conjugators that we find on each such
interval. We then need to understand what happens outside $[L,R]$.

Let $p$ be the rightmost dyadic point of $\partial \Fix(y)\cap [L,R]$. If $y,z \in
\EP_2^>(p,+\infty)$ (or $y,z \in \EP_2^<(p,+\infty)$), then we deal with this case in
Subsection \ref{sec:rescaling-the-circle}. Otherwise, let $q$ be the leftmost point of
$\partial \Fix(y) \cap (R,+\infty)$. If $y(t)=z(t)=t$ on $[p,q]$, then we define $g(t)=t$
on $[p,+\infty)$ and this defines a conjugator for $y$ and $z$ on $[p,+\infty)$ which we
can glue to the previous intervals. Otherwise, we apply the standard conjugacy problem on
the interval $[p,q]$ with final slope $1$ at $q^-$
since the conjugator $g$ has to be the identity
translation on $[R,+\infty)$. If the standard conjugacy problem on $[p,q]$ has no solution,
then $y$ and $z$ cannot be conjugate. Otherwise, if $h$ is the conjugator on $[p,q]$ we
define
\[
g(t):=
\begin{cases}
h(t) & t \in [p,q] \\
t & [q,+\infty)
\end{cases}
\]
which is a well-defined map of $F$, since $g'(q^-)=g'(q^+)=1$, regardless of whether or not $q$ is dyadic.
The map $g$ defines a conjugator for $y$ and $z$ on
$[p,+\infty)$ which we can glue to the previous intervals. A similar argument can be
applied to the left of $L$.

\subsection{Orientation preserving case of the TCP: Mather invariants\label{sec:rescaling-the-circle}}

The procedure outlined in \cite{matucci5} to solve the conjugacy problem in
Bieri-Thompson-Stein-Strebel groups requires various steps which we have studied already:
(i) making $\Fix(y)$ and $\Fix(z)$ coincide (seen in Subsection \ref{ssec:fixed-points})
and (ii) showing that, for a possible initial slope of a conjugator in $F$ (see Remark
\ref{thm:remark-initial-slope}), there exists at most one candidate and we can compute it
through an algorithm (seen in Subsection \ref{ssec:periodicity-boxes}). The next natural
step is to bound the number of integers $\lim_{t \to -\infty}g(t)-t$ representing possible
initial slopes for which we need to build a candidate conjugator.

In order to do this, we will employ ideas to characterize conjugacy seen in
\cite{matucci3}, by taking very large powers of $y$ and $z$ and building a conjugacy
invariant. In \cite{matucci3} a conjugacy class in $F$ has been described by a double coset
$Ay^\infty B$ where $y^\infty$ is an element of Thompson's group $T$ obtained by taking
suitable high powers of $y$ and $A$ and $B$ are two finite cyclic groups (of rotations of
the circle). In the case of the twisted conjugacy problem that we are studying, the Mather
invariant will be essentially defined by a product $A y^\infty B$ where $A \cong B \cong
\mathbb{Z}$.


\emph{Mather invariant construction.} In what follows, we will assume that $y,z \in
\EP_2^>$, to simplify the notation. We can define Mather invariants in the two
neighborhoods of infinity (that is on $\EP_2(-\infty,p)$ and $\EP_2(q,+\infty)$ for some
suitable numbers $p,q$), while solving the conjugacy problem between any two consecutive
dyadic points of $\partial \Fix(y)=\partial \Fix(z)$.

Let $y,z \in \EP_2^>$ and assume that, on the intervals $(-\infty,L]\cup [R,+\infty)$, the
maps $y$ and $z$ coincide and are periodic, for some integers $L \leqslant R$. Let $N \in
\mathbb{N}$ be large enough so that $y^N((y^{-1}(L),L)) \subseteq (R,+\infty)$ and
$z^N((y^{-1}(L),L)) \subseteq (R,+\infty)$.
We look for an orientation preserving homeomorphism $H \in \PL_2(\mathbb{R})$ such that
\begin{enumerate}
\item $H(y^k(L))=k$, for any integer $k$, and
\item $H(y(t))=\lambda(H(t))=H(t)+1$, where $\lambda(t)=t+1$.
\end{enumerate}
To construct $H$, choose any $\PL_2$-homeomorphism $H_0:[y^{-1}(L),L] \to [-1,0]$ with
finitely many breakpoints. Then we extend it to a map $H \in \PL_2(\mathbb{R})$ by defining
 $$
H(t):=H_0(y^{-k}(t))+k \qquad \mbox{if} \; t \in [y^{k-1}(L),y^k(L)] \; \mbox{for some integer $k$}.
 $$
We make a series of remarks.
\begin{itemize}
\item By construction, it is easy to see that $H(y(t))=\lambda(H(t))$ for any real number
    $t$.
\item If we define $\ov{y}:=Hy H^{-1},\ov{z}:=Hz H^{-1}$, we observe that, by
    construction, they both coincide with $\lambda(t)=t+1$ on the intervals
    $(-\infty,1]\cup [N,+\infty)$. It is also clear that $\ov{y}=\lambda$.
\item We notice that $\ov{\lambda}=H\lambda H^{-1} \in \EP_2$. To show this, let $t$ be
    positive sufficiently large so that $y$ is periodic of period $1$ and that all the
    calculations below make sense and define $\widetilde{t}=H_0^{-1}(t- k-1)$:
\[
\ov{\lambda}(t+1)=
H \lambda y^{k+2}H_0^{-1}\lambda^{-k-2}(t+1)=
H\lambda y^{k+2}(\widetilde{t})= H(y^{k+2}(\widetilde{t})+1) =
\]
\[
\lambda^{k'}H_0 y^{-k'}(y^{k+2}(\widetilde{t}+1)) =
\lambda^{k'-1}H_0 y^{-k'+1}(y^{k+1}(\widetilde{t}+1))+1
=\ov{\lambda}(t)+1,
\]
where $k'$ are the jumps that $y$ must make to bring $y^{k+2}(\widetilde{t}+1)$ back to
the domain of $H_0$. A similar argument can be shown for $t$ negative sufficiently large.
\end{itemize}

We define
 $$
C_0:= (-\infty,0)/\mathbb{Z} \qquad C_1:= (N,\infty)/\mathbb{Z}
 $$
and let $p_0:(-\infty,0) \to C_0$ and $p_1:(N,\infty) \to C_1$ be the natural projections.
Then we define the map $\ov{y}^\infty:C_0 \to C_1$ by
 $$
\ov{y}^{\infty}([t]):=[\ov{y}^{N}(t)].
 $$
Similarly we can define the map $\ov{z}^{\infty}$. The maps $\ov{y}^{\infty}$ and
$\ov{z}^{\infty}$ are well-defined and they do not depend on the specific $N$ chosen (the
proof is analogous to the one in Section 3 in \cite{matucci3}). They are called the
\emph{Mather invariants} of $\ov{y}$  and $\ov{z}$ (compare this with the definitions in
Section 3 in \cite{matucci3}).

This induces the equation $\ov{g}\ov{z}^N=\ov{y}^N\ov{g}$ which, following \cite{matucci3},
passes to quotients and becomes
\begin{equation}\label{eq:mather-equivalence}
v_1^k \ov{z}^{\infty} = \ov{y}^{\infty} v_0^{\ell}
\end{equation}
since all the maps $\ov{y},\ov{z},\ov{g}$ are in $\EP_2$ and where $v_1 := p_1 \ov{\lambda}
p_1^{-1}$ is an element of Thompson's group $T_{C_1}$ defined on the circle $C_1$ and
induced by $\ov{\lambda}$ on $C_1$ by passing to quotients via the map $p_1$, $v_0:=p_0
\ov{\lambda} p_0^{-1}$ and where $\ell,k$ are the initial and final slopes of $g$.

The following result shows that the integer solutions of equation
(\ref{eq:mather-equivalence}) correspond to conjugators between $y$ and $z$. The proof is
an extension of the proof of Theorem 4.1 in \cite{matucci3}.

\lemma{\label{thm:mather-iff} Let $y,z \in \EP_2^>$. Then $y$ and $z$ are conjugate through
an element $g \in F$ if and only if there is a pair of integers $k,\ell$ that satisfy
equation (\ref{eq:mather-equivalence}).}

\begin{proof}
Clearly, if $g \in F$ conjugates $y$ to $z$, then equation (\ref{eq:mather-equivalence}) is
satisfied by the calculations above. Conversely, assume that we have a pair $(k,\ell)$ such
that (\ref{eq:mather-equivalence}) is satisfied. We use Theorem
\ref{thm:explicit-conjugator} to find a map $g \in \PL_2(\R)$ which is affine around
$-\infty$, such that $\lim_{x\to -\infty}g(x)-x=\ell$ and that $yg=gz$. By conjugating via
$H$ we see that $\ov{y}\ov{g}=\ov{g}\ov{z}$. If $x$ is positive sufficiently large then
$\ov{y}(x)=\ov{z}(x)=x+1$ so
\[
\ov{g}(x)+1=\ov{y}\ov{g}(x)=\ov{g}\ov{z}(x)=\ov{g}(x+1).
\]
Arguing similarly at $\infty$ we deduce that $\ov{g} \in \EP_2$ and so the equation
$\ov{y}^N\ov{g}=\ov{g}\ov{z}^N$ passes to quotients and becoming $\ov{g}_{\mathrm{ind}}
\ov{z}^{\infty} = \ov{y}^{\infty} v_0^{\ell}$. By using our assumption we see that
$\ov{g}_{\mathrm{ind}} \ov{z}^{\infty} = \ov{y}^{\infty} v_0^{\ell} = v_1^k
\ov{z}^{\infty}$ and by cancellation we obtain $\ov{g}_{\mathrm{ind}} = v_1^k$. By taking
the unique lift of $\ov{g}_{\mathrm{ind}}$ and $v_1^k$ defined on $[N,N+1)$ and passing
through the point $(N,g(N))$, we see that $\ov{g}$ and $\ov{\lambda}^k$ coincide on
$[N,N+1]$ and therefore they coincide on $[N,+\infty)$ since they are both in $\EP_2$.
Thus, $g \in F$ since $\ov{g}(x)=\ov{\lambda}(x)$ around $+\infty$.
\end{proof}

We relabel $t_0:=\ov{z}^{\infty} v_0^{-1} (\ov{z}^{\infty})^{-1}$, $t_1:=v_1$ and and
$t:=\ov{y}^{\infty}(\ov{z}^{\infty})^{-1}$ and we rewrite equation
(\ref{eq:mather-equivalence}) as
\begin{equation}\label{eq:thesis-case}
\numberwithin{equation}{section} t_1^k t_0^\ell=t
\end{equation}
where $t_0,t_1,t \in T_{C_1}$.
To solve equation (\ref{eq:thesis-case}) we will need Lemma 8.4 from \cite{matucci5} which
we restate for the reader's convenience.

\lemmaname{Kassabov-Matucci, \cite{matucci5}}{ \label{thm:matucci5-special-case} Let $p \in
\mathbb{Q}$ and let $\PL_2([p,p+1])$ be the group of piecewise-linear homeomorphisms of the
interval $[p,p+1]$ with finitely many breakpoints which occur at dyadic rational points and
such that all their slopes are powers of $2$. If $t_0,t_1,t \in \PL_2([p,p+1])$, there is
an algorithm which outputs one of the following two mutually exclusive cases in finite
time:
\begin{enumerate}
\item Equation (\ref{eq:thesis-case}) has at most one solution and we compute a pair
    $(k,\ell)$ such that, if equation (\ref{eq:thesis-case}) is solvable, then $(k,\ell)$
    must be its unique solution.
%
\item Equation (\ref{eq:thesis-case}) has infinitely many solutions which are given by
    the sequence of pairs $(k_j,\ell_j)$ where $k_j = a_1 j +b_1$ and $\ell_j=a_2 j +
    b_2$ for any $j \in \mathbb{Z}$ and for some integers $a_1,a_2,b_1,b_2$ which we can
    compute.
\end{enumerate}
}

Lemma \ref{thm:matucci5-special-case} gives a solution for equation (\ref{eq:thesis-case})
in the case that $t_0,t_1,t$ live in a copy of Thompson's group $\PL_2([p,p+1])$ of
functions over an interval. However, equation (\ref{eq:thesis-case}) needs to be solved in
a copy of Thompson's group $T$ of functions over a circle, so we will need to adapt Lemma
\ref{thm:matucci5-special-case} to our needs.

%
%
%
%

\lemma{ \label{thm:thesis-case} Let $T$ be Thompson's group $\PL_2(S^1)$ and let $t_0,t_1,t
\in T$. Then there is an algorithm which outputs one of the following two mutually
exclusive cases in finite time:
\begin{enumerate}
\item Equation (\ref{eq:thesis-case}) has at most finitely many solutions and we compute
    a finite set $S$ such that, if $(k,\ell)$ is a solution of equation
    (\ref{eq:thesis-case}), then $\ell \in S$.
\item Equation (\ref{eq:thesis-case}) has infinitely  infinitely many solutions and we
    compute a sequence of solutions $(k_j,\ell_j)$ where $k_j = a_1 j +b_1$ and
    $\ell_j=a_2 j + b_2$ for any $j \in \mathbb{Z}$ and for some integers
    $a_1,a_2,b_1,b_2$.

\end{enumerate}
}

\begin{proof}
For a map $h \in T$, we denote by $\mathrm{Per}(h)$ the set of all periodic points of $h$.
Obviously, $\Fix(h) \subseteq \mathrm{Per}(h)$. By a result of Ghys and Sergiescu
\cite{GhysSergiescu} every element of $T$ has at least one periodic
point. For $i=0,1$, we find a $q_i \in \mathrm{Per}(t_i)$ be a point of period $d_i$. If
$d=\mathrm{lcm}(d_0,d_1)$, then both $t_0^d,t_1^d$ have fixed points and therefore
$\mathrm{Per}(t_i^d)=\Fix(t_i^d)$.

Using the division algorithm we write $k=k'd+r$ and $\ell=\ell'd+s$ with $0 \leqslant r,s
<d$ so that equation (\ref{eq:thesis-case}) becomes

\begin{equation}\label{eq:family-of-equations}
\numberwithin{equation}{section} (t_1^d)^{k'}(t_0^d)^{\ell' }=t_1^{-r}t t_0^{-s}.
\end{equation}
By considering all possibilities for $0 \leqslant r,s <d$, equation
(\ref{eq:family-of-equations}) can be regarded as a family of $d^2$ equations in $T$.
Equation (\ref{eq:thesis-case}) is solvable if and only if at least one of the $d^2$
equations (\ref{eq:family-of-equations}) is solvable.

Up to renaming $t_0^d$ with $t_0$, $t_1^d$ with $t_1$ and $t_1^{-r}t t_0^{-s}$ with $t$, we
observe that each of the equations (\ref{eq:family-of-equations}) has the same form of
equation (\ref{eq:thesis-case}), therefore we have reduced ourselves to study equation
(\ref{eq:thesis-case}) with the extra assumption that both $t_0$ and $t_1$ have fixed
points. We compute the full fixed point sets of $t_0$ and $t_1$. We now break the proof
into two cases.

\noindent \emph{Case 1: There is a point $p \in \partial \Fix(t_1)$ such that $p \not \in
\Fix(t_0)$.} Rewriting equation (\ref{eq:thesis-case}) and applying it to $p$, we get
\begin{equation}\label{eq:thesis-cases-rewrite}
\numberwithin{equation}{section} t_0^{-\ell}(p)= t^{-1}(p).
\end{equation}
Since $p \not \in \Fix(t_0)$ and $t_0$ is orientation preserving, then there exists at most
one number $\ell$ satisfying equation (\ref{eq:thesis-cases-rewrite}) by Lemma
\ref{thm:standard-folklore}.

\noindent \emph{Case 2: There is a point $p \in \partial \Fix(t_1) \cap \Fix(t_0)$.} If $p
\not \in \Fix(t)$, by particularizing at $p$ we see that equation (\ref{eq:thesis-case}) is
not solvable for any pair $(k,\ell)$. Otherwise, $p \in \Fix(t)$ and we can cut the unit
circle open at $p \in \mathbb{Q}/\mathbb{Z}$ and regard $t_0,t_1,t$ as elements of
$\PL_2([p,p+1])$. We can now finish the proof by applying Lemma
\ref{thm:matucci5-special-case}.
\end{proof}

\remark{The proof of Lemma \ref{thm:thesis-case} shows how to locate the pairs $(k,\ell)$.
We need to find all periodic orbits and their periods and this can be effectively achieved
by computing the Brin-Salazar revealing pairs of the tree pair diagrams of $T$, using the
Brin-Salazar technology to compute neutral leaves and thus deducing the size of periodic
orbits (see, for example, Section 4 in \cite{matucci8}).
}

\remark{ \label{thm:counting-admissible-slopes} To sum up this subsection, Lemma
\ref{thm:mather-iff} shows that $y$ and $z$ are conjugate via an element of $F$ if and only
if equation (\ref{eq:mather-equivalence}) is solvable for some integers $k,\ell$. Lemma
\ref{thm:thesis-case} shows how to narrow down the number of pairs $(k,\ell)$ that we need
to test. There are two possible cases:
\begin{enumerate}
\item[(i)] In the first case of Lemma \ref{thm:thesis-case} we are given a finite set $S$
    of initial slopes to test. We can use Theorem \ref{thm:explicit-conjugator}(ii) for
    each of the finitely many initial slopes in the set $S$. There is a conjugator if and
    only if one of the applications of Theorem \ref{thm:explicit-conjugator}(ii) returns
    a positive answer. If there is a conjugator, it can be built using Theorem
    \ref{thm:explicit-conjugator}(i).
\item[(ii)] In the second case of Lemma \ref{thm:thesis-case} there are infinitely many
    possible pairs $(k,\ell)$ (and we can construct explicitly an infinite family) and
    all of them correspond to a conjugator between $y$ and $z$. We can apply Theorem
    \ref{thm:explicit-conjugator}(i) on a specific pair $(k,\ell)$ of our choice to find
    an explicit conjugator between $y$ and $z$
\end{enumerate}
Hence, in every case we can find at least one conjugator, if it exists.}

\remark{We observe that the construction of the Mather invariant can be carried out even
when $y$ and $z$ are elements of $\EP_2^>(p,+\infty)$ or of $\EP_2^>(-\infty,p)$ for any
rational number $p$. All the results of the current subsection can still be recovered. For
this reason, in the following we will refer to the Mather invariant regardless of the
ambient set where it will be built.}

\subsection{Orientation reversing case of the TCP\label{sec:special-case}}

We now study the orientation reversing case of TCP, that is, we want to solve the equation
\begin{equation}\label{eq:conj}\numberwithin{equation}{section}
z=g^{-1}yg,
\end{equation}
where $y,z \in \REP \setminus \{\id\}$ and $g \in F$. The general idea that we will follow
is to square the equation and attempt to solve
\[
z^2 = g^{-1} y^2 g
\]
so that $y^2,z^2 \in \EP_2$ and we can appeal to the results of the previous subsections.

Since $y,z$ are strictly decreasing and approach $\mp \infty$ when $t \to \pm \infty$ then
both $y$ and $z$ have exactly one fixed point each. Moreover, all possible $g$'s fulfilling
equation~(\ref{eq:conj}) must also satisfy $g(\Fix(z))= \Fix(gzg^{-1})=\Fix(y)$. By
Proposition~\ref{thm:rationals-coincide}(ii), one can algorithmically decide whether or not
there is $g \in F$ mapping the point $\Fix(z)$ to the point $\Fix(y)$. If there is no such
$g$, then equation~(\ref{eq:conj}) has no solution and we are done. Otherwise, compute such
a $g\in F$ and, after replacing $z$ by $gzg^{-1}$, we can assume that $\Fix(y)=\Fix(z)=\{
p\}$, for some $p\in \mathbb{Q}$.

We start with a special case and then move on to consider all orientation reversing maps.

\proposition{\label{thm:orientation-reversing-special}
Let $y,z \in \mathcal{R}\cdot \EP_2$ be such that $y^2=z^2=\id$ and $y(p)=z(p)=p$, for some
$p\in \mathbb{Q}$. Then $y$ and $z$ are conjugate by an element of $F$ if and only if there
exists
$u\in \mathbb{Z}$ such that
$y^{-1}z(t)=t+u$ for $t$ positive sufficiently large. }
\begin{proof}
The forward direction follows from a straightforward check of the behavior of $y$ and $z$
at neighborhoods of $\pm \infty$. For the converse, define the following map
 $$
g(t):=
\begin{cases}
t & \; \mbox{if} \; \; \; \; \; t \in (-\infty,p] \\ y^{-1}z(t) & \; \mbox{if} \; \; \; \; \; t \in [p,+\infty).
\end{cases}
 $$
If $t\leqslant p$, then
 $$
g(t)=t=y^{-2}z^2(t)=y^{-1}(y^{-1}z)z(t)=y^{-1}gz(t)
 $$
since $y^2=z^2=\id$ and $z(t)\geqslant p$. On the other hand, if $t\geqslant p$, then
 $$
g(t)=y^{-1}z(t)= y^{-1}gz(t)
 $$
since $z(t) \leqslant p$. So $y$ and $z$ are conjugate to each other by the element $g\in
\EP_2$. The final step is to observe that $g$ is, in fact, in $F$ because $g(t)=t$, for $t$
negative sufficiently large, and $g(t)=t+u$ by construction, for $t$ positive sufficiently
large.
\end{proof}

We quickly recall and extend an argument from \cite{matucci5} to reduce the number of
candidate conjugators to test. The trick is to reduce the number of initial slopes that we
need to test.

\lemma{ \label{thm:solving-reverse-squared} Let $\ov{y}, \ov{z} \in \REP^<(p,+\infty)$ and
$g \in F(p,+\infty)$ and consider the equation
\begin{equation}\label{eq:standard-conjugacy-equation}
\numberwithin{equation}{section} \ov{z} = x^{-1} \ov{y} x.
\end{equation}
Then $x=g$ is a solution of (\ref{eq:standard-conjugacy-equation}) if and only if there
exists an integer $n$ such that $x=\ov{y}^{2n}g \in \EP_2(p,+\infty)$ is the unique
solution of equation (\ref{eq:standard-conjugacy-equation}) such that $(\ov{y}^{2n}g)'(p^+)
\in [(y^2)'(p^+),1]$. }

\begin{proof}
This follows immediately by noticing that equation (\ref{eq:standard-conjugacy-equation})
is equivalent to
 $$
\ov{z} = (\ov{y}^{2n}x)^{-1}\ov{y} (\ov{y}^{2n} x).
 $$
To show uniqueness, we observe that in Subsection \ref{sec:restatement-TCP} we noticed that
a solution of equation (\ref{eq:standard-conjugacy-equation}) is also a solution of the
squared equation
\begin{equation}\label{eq:squared-conjugacy-equation-first}
\numberwithin{equation}{section} \ov{z}^2 = g^{-1} \ov{y}^2 g.
\end{equation}
Uniqueness follows from Theorem \ref{thm:explicit-conjugator} applied to the squared
equation (\ref{eq:squared-conjugacy-equation-first}).
\end{proof}

\theorem{\label{thm:orientation-reversing-general} Let $y,z \in \mathcal{R}\cdot \EP_2$ be
such that $y(p)=z(p)=p$, for some $p\in \mathbb{Q}$. We can decide whether or not $y$ and
$z$ are conjugate by an element of $F$. If there exists a conjugator, we can construct
one.}
\begin{proof}
If $y^2=z^2=\id$, then we are done by Proposition \ref{thm:orientation-reversing-special}.
Moreover, if $y$ and $z$ are conjugate via an element of $F$, it is immediate that
$y^{-1}z(t)=t+u$, for some integer $u$ and for any $t$ positive sufficiently large (as
observed in the proof of Proposition \ref{thm:orientation-reversing-special}). Thus we can
assume that $y^{-1}z$ is a translation, for $t$ positive sufficiently large.

Assume now $y^2\ne \id \ne z^2$. We can appeal to Proposition
\ref{thm:identify-fixed-points} and assume that $\Fix(y^2)=\Fix(z^2)$, up to suitable
conjugation. Moreover, if there exists a conjugator between $y$ and $z$, then it must fix
$\Fix(y)=\Fix(z)=\{p\}$ and so $\{p\} \subseteq \Fix(y^2)=\Fix(z^2) \subseteq \Fix(g)$.

Let $L<R$ be two suitable integers so that $y^2$ and $z^2$ coincide and are periodic on the
set $(-\infty,L] \cup [R,+\infty)$. If either $L$ or $R$ does not exist, then $y$ and $z$
cannot be conjugate. We can apply the techniques from \cite{matucci5} on any two
consecutive dyadic rational points $p_1,p_2$ of $\partial \Fix(y^2) \cap [L,R]$ where
$y^2|_{[p_1,p_2]} \ne \id|_{[p_1,p_2]}$ and $z^2|_{[p_1,p_2]} \ne \id|_{[p_1,p_2]}$ and
find (if they exist) all the finitely conjugators between $y^2|_{[p_1,p_2]}$ and
$z^2|_{[p_1,p_2]}$ with initial slopes within $(y^2)'(p^+)$ and $(y^{-2})'(p^+)$. Similarly
we can do on $[a,+\infty)$ where $a$ is the rightmost dyadic rational point of $\partial
\Fix(y^2) \cap [L,R]$ by applying Lemma \ref{thm:solving-reverse-squared} in the case that
$y^2$ and $z^2$ have no fixed points on $[R,+\infty)$ (to reduce the number of initial
slopes on which we can apply Theorem \ref{thm:explicit-conjugator}) or using the argument
at the end of Subsection \ref{ssec:reducing-to-squares} in case $y^2$ and $z^2$ have fixed
points on $[R,+\infty)$.

Thus in all cases, up to using the same trick of Lemma \ref{thm:solving-reverse-squared} to
reduce the slopes to test, we apply Theorem \ref{thm:explicit-conjugator} (or its bounded
version from \cite{matucci5}) to build finitely many functions between any two consecutive
dyadic rational points $p_1,p_2$ of $\partial \Fix(y^2)$ (respectively, on an interval of
the type $[p_1,+\infty)$) and such that $y^2 \ne \id$ on $[p_1,p_2]$ (respectively, on an
interval of the type $[p_1,+\infty)$).

We now test all these functions as conjugators between $y$ and $z$ in the respective
intervals. If there is an interval such that none of these functions conjugates $y$ and
$z$, then $y$ and $z$ cannot be conjugate via an element of $F$. Otherwise, on each such
interval $U_s$ we fix a conjugator $g_s$ between $y$ and $z$.

Now we will carefully glue all these conjugators with the function that we have built in
Proposition \ref{thm:orientation-reversing-special}. Assume that $(p,+\infty) \setminus
\Fix(y^2)$ is a disjoint union of ordered intervals $I_i=(a_i,b_i)$ so that $a_i < a_j$, if
$i<j$. Similarly, assume that $(-\infty,p) \setminus \Fix(y^2)$ is a disjoint union of
ordered intervals $J_i=(d_i,c_i)$ such that $c_i>c_j$, if $i<j$.
$$
g(t):=
\begin{cases}
t & \; \mbox{if} \; \; \; \; \; t=p \text{ or } t \in \Fix(y^2) \cap (-\infty,p) \\
y^{-1}z(t) & \; \mbox{if} \; \; \; \; \; t \in \Fix(y^2) \cap (p,+\infty) \\
g_s(t) & \; \mbox{if} \; \; \; \; \; t \in U_s
\end{cases}
$$
Since $y$ acts on $\R$ in an order reversing way, it is immediate to verify that
$y(a_i)=c_i=z(a_i)$, $y(c_i)=a_i=y(c_i)$, $y(b_i)=d_i=z(b_i)$ and $y(d_i)=b_i=z(d_i)$ and
therefore the map $g$ is in $F$. It is straightforward to observe that this map is
continuous and in $F$ and that it is a conjugator, by construction. For example, since
$z([c_{i+1},d_i])=[b_i,a_{i+1}]$ and $y^2=z^2=\id$ on $[c_{i+1},d_i]$ then it is clear that
$$
g(t)=t=y^{-2}z^2(t)=y^{-1}(y^{-1}z)z(t)=y^{-1}gz(t)
$$
for any $t \in [c_{i+1},d_i]$.
\end{proof}

\subsection{Solution of the TCP\label{ssec:solution-TCP}}

We are now ready to prove Theorem \ref{thm:TCP-solvable}.

\medskip
\noindent \textbf{Theorem \ref{thm:TCP-solvable}.} \emph{Thompson's group $F$ has solvable
twisted conjugacy problem.}

\begin{proof}
Given $y,z \in F$ and $\varphi \in \Aut(F)$, we need to establish whether or not there is a
$g \in F$ such that
\begin{equation}\label{eq:recall-TCP}
\numberwithin{equation}{section} z = g^{-1} y \varphi(g).
\end{equation}
In Subsection \ref{sec:restatement-TCP} we have shown that equation (\ref{eq:recall-TCP})
is equivalent to the equation
\begin{equation}\label{eq:original-conjugacy-equation}
\numberwithin{equation}{section} \ov{z} = g^{-1} \ov{y} g,
\end{equation}
for $\ov{y},\ov{z} \in \WTEP_2$ and $g \in F$. We describe a procedure to wrap up all work
of the previous subsections:

\begin{enumerate}
\item[(1)] If one of $\ov{y}$ and $\ov{z}$ belongs to $\EP_2$ and the other in $\REP$,
    then equation (\ref{eq:original-conjugacy-equation}) has no solution for $g\in F$,
    since conjugation does not change the orientation of a function.
\item[(2)] If both $\ov{y},\ov{z} \in \EP_2$, then we apply the results of Subsections
    \ref{ssec:periodicity-boxes} through \ref{sec:rescaling-the-circle} to solve equation
    (\ref{eq:original-conjugacy-equation}).
\item[(3)] If $\ov{y},\ov{z} \in \REP$, then we apply Theorem
    \ref{thm:orientation-reversing-general} to solve equation
    (\ref{eq:original-conjugacy-equation}).
\end{enumerate}
This ends the proof of Theorem \ref{thm:TCP-solvable}.
\end{proof}

\section{Extensions of $F$ with unsolvable conjugacy problem \label{sec:CP-extensions}}

In this section we recall the necessary tools from~\cite{bomave2} in order to construct
extensions of Thompson's group $F$ with unsolvable conjugacy problem (proving
Theorem~\ref{thm:CP-extension-unsolvable}).

As explained in the introduction, Bogopolski, Martino and Ventura give a criterion to study
the conjugacy problem in extensions of groups (see Theorem~\ref{thm:bomave-extensions}).
Applying it to the case we are interested in, let $F$ be Thompson's group, let $H$ be any
torsion-free hyperbolic group (for example, a finitely generated free group), and consider
an algorithmic short exact sequence
\begin{equation} \label{eq:exact-sequence}
\numberwithin{equation}{section} 1 \longrightarrow F \overset{\alpha}{\longrightarrow} G
\overset{\beta}{\longrightarrow} H \longrightarrow 1.
\end{equation}
We can then consider the \emph{action subgroup} of the sequence, $A_G =\{\varphi_g \mid g
\in G\} \leqslant \Aut(F)$, and Theorem~\ref{coro} tells us that $G$ has solvable conjugacy
problem if and only if $A_G\leqslant \Aut(F)$ is orbit decidable. In the present section we
will find orbit undecidable subgroups of $\Aut(F)$ and so, extensions of Thompson's group
$F$ with unsolvable conjugacy problem.

A good source of orbit undecidable subgroups in $\Aut(F)$ comes from the presence of
$F_2\times F_2$ via Theorem 7.4 from~\cite{bomave2}:

\theoremname{Bogopolski-Martino-Ventura, \cite{bomave2}}{\label{thm:bomave-unsolvable} Let
$F$ be a finitely generated group such that $F_2 \times F_2$ embeds in $\Aut(F)$ in such a
way that the image $B$ intersects trivially with $\mathrm{Stab}^*(v)$ for some $v \in F$,
where
 $$
\mathrm{Stab}^*(v)=\{ \theta \in \Aut(F) \mid \theta(v) \text{ is conjugate to } v \text{ in } F\}.
 $$
Then $\Aut(F)$ contains an orbit undecidable subgroup.}

Let us first find a copy of $F_2\times F_2$ inside $\Aut(F)$ and then deal with the
technical condition about avoiding the stabilizer.

We can define two maps $\varphi_{-\infty},\, \varphi_{\infty} \colon \EP_2 \to
T=\PL_2(S^1)$ in the following way: given $f\in \EP_2$ we find a negative sufficiently
large integer $L$ so that $f$ is periodic in $(-\infty,\, L]$; then we pass $f|_{(L-1,\,
L]}$ to the quotient modulo $\mathbb{Z}$ to obtain an element from $T$ defined to be the
image of $f$ by $\varphi_{-\infty}$. The map $\varphi_{+\infty}$ is defined similarly but
looking at a neighborhood of $+\infty$.

The maps $\varphi_{-\infty}$ and $\varphi_{+\infty}$ are clearly well-defined homomorphisms
from $\EP_2$ to $T$. Note also that, for $f_1, f_2\in \EP_2$ and $k\in \mathbb{Z}$, if
$f_1$ and $f_2+k$ agree for $t$ negative (resp. positive) sufficiently large, then
$\varphi_{-\infty}(f_1)=\varphi_{-\infty}(f_2)$ (resp.
$\varphi_{+\infty}(f_1)=\varphi_{+\infty}(f_2)$).

We begin by showing that both $\varphi_{-\infty}$ and $\varphi_{+\infty}$ are surjective.

\lemma{\label{thm:embed-F_2} For every $a\in T$ and every dyadic rational $p$, there
exist preimages of $a$ by $\varphi_{-\infty}$ and $\varphi_{+\infty}$, respectively inside
$\EP_2(-\infty,p) \leqslant \EP_2$ and $\EP_2(p,+\infty) \leqslant \EP_2$.}

\begin{proof}
We show the result for the case $\EP_2(p,+\infty)$ (the other case is completely
analogous). Let $a\in T$ and choose $\wt{a}\in \EP_2$ to be any standard periodic lift of
$a$ conveniently translated up so that $p<\wt{a}(p+1)$. By Proposition~\ref{thm:rationals-coincide},
we can construct $g\in F$ such that $g(p)=p$ and $g(p+1)=\wt{a}(p+1)$. Finally, consider
 $$
\wh{a}(t)= \begin{cases} t & t\leqslant p \\ g(t) & p\leqslant t\leqslant p+1 \\ \wt{a}(t) &
p+1\leqslant t,
\end{cases}
 $$
which is clearly an element of $\EP_2(p,+\infty)$ such that $\varphi_{+\infty}(\wh{a})
=\varphi_{+\infty}(\wt{a})=a$.
\end{proof}

The following Corollary is the key observation of the current subsection.

\corollary{\label{thm:F_2 x F_2 in EP_2} The automorphism group of Thompson's group
$F=\PL_2(I)$ contains a copy of the direct product of two free groups, $F_2 \times F_2
\leqslant \EP_2 \leqslant \Aut^+(F)$.}

\begin{proof}
It is well known that Thompson's group $T=\PL_2(S^1)$ contains a copy of $F_2$, the free
group on two generators, say generated by $a,b\in T$. Apply Lemma~\ref{thm:embed-F_2} to
obtain preimages of $a$ and $b$ by $\varphi_{-\infty}$, say $\wh{a}_-, \wh{b}_- \in
\EP_2(-\infty ,0)$, and preimages of $a$ and $b$ by $\varphi_{+\infty}$, say $\wh{a}_+,
\wh{b}_+ \in \EP_2(0, +\infty)$. Since $\varphi_{-\infty}$ and $\varphi_{+\infty}$ are
homomorphisms, we have again $\langle \wh{a}_-, \wh{b}_-\rangle \simeq F_2 \simeq \langle
\wh{a}_+, \wh{b}_+\rangle$. And, on the other hand, by disjointness of supports, they
commute to each other and so $F_2 \times F_2 \simeq \langle \wh{a}_-, \wh{b}_-, \wh{a}_+,
\wh{b}_+\rangle \leqslant \EP_2 \simeq \Aut^+(F)$.
\end{proof}

We are finally ready to prove Theorem \ref{thm:CP-extension-unsolvable}.

\medskip
\noindent \textbf{Theorem \ref{thm:CP-extension-unsolvable}.} \emph{There are extensions of
Thompson's group $F$ by finitely generated free groups, with unsolvable conjugacy problem.}

\begin{proof}
We need to redo the proof of Corollary~\ref{thm:F_2 x F_2 in EP_2} in an algorithmic
fashion and choosing our copy of $F_2\times F_2$ inside $\Aut^+(F)$ carefully enough so
that it satisfies the technical condition in Theorem~\ref{thm:bomave-unsolvable}.

Let $\Theta$ be the map obtained
by repeating periodically the map $\theta$ defined in Subsection \ref{sec:thompson and autos} inside 
each square $[k,k+1]^2$, for any integer $k$.
Let $\alpha(t):= \Theta^2(t) \pmod{1} \in T$ and $\beta(t):=\Theta^2(t)+\frac{1}{2} \pmod{1} \in T$.
By using the ping-pong lemma it is straightforward to verify that $\alpha$ and $\beta$
generate a copy of $F_2$ inside $T$. Now take $a=\alpha^2$,
$b=\beta^2$, $c=\alpha \beta \alpha^{-1}$ and $d=\beta \alpha \beta^{-1}$, which generate a
copy of the free group of rank four, $F_4 \simeq \langle a,b,c,d\rangle \leqslant T$.

Using Lemma~\ref{thm:embed-F_2}, we can find preimages of $a,b\in T$ by
$\varphi_{-\infty}$, denoted by $\wh{a},\, \wh{b}\in \EP_2(-\infty,\, 0) \leqslant \EP_2$,
and preimages of $c,d\in T$ by $\varphi_{+\infty}$, denoted by $\wh{c},\, \wh{d} \in
\EP_2(0,\, +\infty) \leqslant \EP_2$. Since $\langle a,b\rangle \cong F_2 \cong \langle
c,d\rangle$ and $\varphi_{-\infty}$ and $\varphi_{+\infty}$ are both group homomorphisms,
we get $\langle \wh{a},\wh{b}\rangle \cong F_2 \cong \langle \wh{c},\wh{d}\rangle$.
Moreover, the disjointness of supports gives us that $F_2\times F_2 \cong \langle
\wh{a},\wh{b},\wh{c}, \wh{d}\rangle \leqslant \EP_2$; this is the copy $B$ of $F_2\times
F_2$ inside $\EP_2$ (though as positive automorphisms of $F$ via Brin's Theorem) ready to
apply Theorem~\ref{thm:bomave-unsolvable}. Additionally, note that, by construction,
$\varphi_{-\infty}(\wh{a})=a$, $\varphi_{-\infty}(\wh{b})=b$, $\varphi_{+\infty}(\wh{c})=c$
and $\varphi_{+\infty}(\wh{d})=d$ but, at the same time,
$\varphi_{+\infty}(\wh{a})=\varphi_{+\infty}(\wh{b})=
\varphi_{-\infty}(\wh{c})=\varphi_{-\infty}(\wh{d})=1_T$.

Let now $v\in F$ be the map $v(t)=t+1$, for all $t \in \mathbb{R}$. We will show that
$B\cap \mathrm{Stab}^*(v)=\{ \id \}$. Let $\tau \in B \cap \mathrm{Stab}^*(v)$. On one
hand, $\tau \in B$ and so $\tau(0)=0$ and $\tau = w_1(\wh{a},\wh{b})w_2(\wh{c},\wh{d})$ for
some unique reduced words $w_1(\wh{a},\wh{b}) \in \langle \wh{a},\wh{b} \rangle$ and
$w_2(\wh{c},\wh{d}) \in \langle \wh{c},\wh{d} \rangle$. On the other hand, $\tau \in
\mathrm{Stab}^*(v)$ and so $\tau^{-1} v\tau = g^{-1}v g$ for some $g\in F$, which implies
that $\tau g^{-1}$ commutes with $v$ in $\EP_2$. By definition of $v$, the map $\tau
g^{-1}$ is periodic of period 1 on the entire real line, thus $\varphi_{-\infty}(\tau
g^{-1})=\varphi_{+\infty}(\tau g^{-1})$ in $T$. On the other hand, since $g\in F$, there
exist integers $m_-$ and $m_+$ such that, for negative sufficiently large $t$, $\tau
g^{-1}(t)=\tau(t-m_-)=\tau(t)-m_-$, and for positive sufficiently large $t$, $\tau
g^{-1}(t)=\tau(t-m_+)=\tau(t)-m_+$. Modding out these two equations by $\mathbb{Z}$ around
$\pm \infty$, we get
 $$
\varphi_{-\infty}(\tau g^{-1})=\varphi_{-\infty}(\tau )=\varphi_{-\infty} (w_1(\wh{a},\wh{b})
w_2(\wh{c},\wh{d}))=
 $$
 $$
=\varphi_{-\infty} (w_1(\wh{a},\wh{b})) \varphi_{-\infty}(w_2(\wh{c},\wh{d}))=w_1(a,b);
 $$
similarly, $\varphi_{+\infty}(\tau g^{-1})=w_2(c,d)$. Hence,
 $$
w_1(a,b)=\varphi_{-\infty}(\tau g^{-1})=\varphi_{+\infty}(\tau g^{-1})=w_2(c,d),
 $$
an equation holding in $\langle a,b,c,d\rangle \leqslant T$. Since this is a free group on
$\{ a,b,c,d\}$, we deduce that $w_1(a,b)$ and $w_2(c,d)$ are the trivial words and
therefore $\tau =\id$.

Having shown that $B\cap \mathrm{Stab}^*(v)=\{\id\}$, an application of
Theorem~\ref{thm:bomave-unsolvable} gives us orbit undecidable subgroups of $\Aut^+(F)$,
and Theorem~\ref{coro} concludes the proof.
\end{proof}

\remark{The element $v$ chosen in the previous proof is actually $x_0$, the first generator
of the standard finite presentation defined in Subsection \ref{sec:thompson and autos}.}

\section{The orbit decidability problem for $F$ \label{sec:ODP-solvable}}

In this section we study the orbit decidability problem for $\Aut(F)$ and $\Aut_+(F)$. We
study two different cases and use techniques which are ``dual'' to those of Section
\ref{sec:twisted-problem}. As a consequence, provided that one knows the solvability of a
certain decision problem, we can build nontrivial extensions of $F$ with solvable conjugacy
problem.

By using Theorem \ref{thm:brin-thm} and following computations similar to those in
Subsection \ref{sec:restatement-TCP}, the orbit decidability problem for $\Aut(F)$ can be
restated as the following one: given $y,z \in F$ decide whether or not there exists a $g
\in \EP_2$ such that either
\begin{enumerate}
\item[(i)] $g^{-1}yg=z$, or
\item[(ii)] $g^{-1}(\mathcal{R}y\mathcal{R})g=z$.
\end{enumerate}
Notice that the first equation corresponds to orbit decidability for $\Aut_+(F)$. Up to
renaming $\mathcal{R}y\mathcal{R}$ by $y$, both (i) and (ii) can be regarded as an instance
of (i).

\subsection{Orbit decidability problem: fixed points\label{ssec:ODP-fixed-points}}

It is immediate to adapt Lemma \ref{thm:identical-at-infinity} to this setting, noticing
that if $y \sim_{\EP_2} z$ then $y$ and $z$ coincide around $\pm \infty$.



\remark{ \label{thm:remark-fixed-points-coincide} Since $y,z \in F$ have only finitely many
intervals of fixed points, we can use the results of Subsection \ref{ssec:fixed-points} and
assume that $\Fix(y)=\Fix(z)$, up to conjugating by a $g \in F$. It can be shown that if
there is no $g \in F$ such that $\Fix(y)=g(\Fix(z))$, then there is no $h \in \EP_2$ such
that $\Fix(y)=h(\Fix(z))$.}


\lemma{ \label{thm:ODP-with-fixed-points} Let $y,z \in F$ such that $\Fix(y)=\Fix(z) \ne
\emptyset$. It is decidable to determine whether or not there is a $g \in \EP_2$ such that
$g^{-1}yg=z$. }

\begin{proof}
If $g \in \EP_2$ conjugates $y$ to $z$, then it must fix $\Fix(z)$ point wise.
For any two consecutive points $p_1,p_2$ of $\partial \Fix(z)$ we can use the techniques in
\cite{matucci5} to decide whether or not there is a $h_{p_1,p_2} \in \PL_2([p_1,p_2])$
conjugating $y|_{[p_1,p_2]}$ to $z|_{[p_1,p_2]}$.

Let $R=\max \Fix(z)$. If $R=+\infty$, then there exists a rational number $p$ such that
$y=z=\id$ on $[p,+\infty)$ and so we can choose $g \in \EP_2(R,+\infty)$ to be $g=\id$ to
conjugate $y$ to $z$. Assume now that $R<+\infty$.

By using the same idea seen in Subsection \ref{ssec:solution-TCP} and rewriting the
equation $z=g^{-1}yg=(y^n g)^{-1} y (y^n g)$ we restrict to looking for candidate
conjugators with slopes at $R^+$ inside $[y'(R^+),1]$. For any power $2^\alpha$ within
$[y'(R^+),1]$, we apply Theorem \ref{thm:explicit-conjugator}(ii) to build the unique
conjugator $g \in \PL_2(R,+\infty)$ such that $g'(R^+) = 2^\alpha$. We find a finite number
of conjugators $g_1,\ldots,g_s \in \PL_2(R,+\infty)$. Notice: by Theorem
\ref{thm:explicit-conjugator}(ii) every $g_i$ conjugates $y$ to $z$, but it may not be true
that $g_i \in \EP_2(R,+\infty)$.

There exists a positive sufficiently large number $M$ such that, for any $t \geqslant M$,
we have $y(t)=t+k=z(t)$ and that for any $i=1,\ldots,s$ and any $t \geqslant M$, we have:
 $$
g_i(t)+k=yg_i(t)=g_iz(t)=g_i(t+k),
 $$
so that every $g_i$ is periodic of period $k$ on $[M,+\infty)$. To finish the proof, we
only need to check if any of the $g_i$'s is in $\EP_2(R,+\infty)$. To do so, we check if
$g_i(t+1)=g_i(t)+1$ on the interval $[M,M+k]$. If any of them is indeed periodic of period
$1$, then we have found a valid conjugator, otherwise $y$ and $z$ are not conjugate.
\end{proof}

\subsection{
\label{ssec:ODP-Mather}Orbit decidability problem: Mather invariants }

We assume that $y,z \in F^>$ and that there exist two integers $L < R$ such that
$y(t)=z(t)=t+a$ for $t \leqslant L$ and $y(t)=z(t)=t+b$ for $t \geqslant R$, for suitable
integers $a,b \geqslant 1$. Up to conjugation by a suitable $g \in F$, we can assume that
$L=0$ and $R=1$. Define the two circles
 $$
C_0:= (-\infty,0)/a\mathbb{Z} \qquad C_1:= (1,\infty)/b\mathbb{Z}
 $$
and let $p_0:(-\infty,0) \to C_0$ and $p_1:(1,\infty) \to C_1$ be the natural projections.
As was done before, let $N$ be a positive integer large enough so that $y^N(-a,0) \subseteq
(1,+\infty)$ and define the map $y^\infty:C_0 \to C_1$ by
 $$
y^{\infty}([t]):=[\ov{y}^{N}(t)].
 $$
Similarly we define $z^{\infty}$ and call them the \emph{Mather invariants} for $y$ and
$z$. Arguing as in Subsection \ref{sec:rescaling-the-circle} we see that, if $g^{-1}yg = z$
for $g \in \EP_2$, then
\begin{equation}\label{eq:ODP-mather}
\numberwithin{equation}{section} v_1 z^{\infty}=y^{\infty} v_0
\end{equation}
where $v_i$ is an element of Thompson's group $T_{C_i}$ induced by $g$ on $C_i$, for
$i=0,1$, and such that $v_i(t+1)=v_i(t)+1$.

Recall that a group $G$ has solvable \emph{$k$-simultaneous conjugacy problem} ($k$-CP) if,
for any two $k$-tuples $(y_1, \ldots, y_k)$, $(z_1, \ldots, z_k)$ of elements of $G$, it is
decidable to say whether or not there is a $g\in G$ so that $g^{-1} y_i g=z_i$, for all
$i=1,\ldots, k$. Kassabov and the second author ~\cite{matucci5} show that Thompson's group
$F$ has solvable $k$-CP.

\conjecture{\label{conj-T} Thompson's group $T$ has solvable $k$-CP. }

This conjecture is believed to be true, and partial results have been obtained by Bleak,
Kassabov and the second author in Chapter~7 of the second author's
thesis~\cite{matuccithesis}; it is work in progress to complete this investigation.

\lemma{ \label{thm:ODP-without-fixed} Let $y,z \in F^>$. If the $2$-simultaneous conjugacy
problem is solvable in Thompson's group $T$, then it is decidable to determine whether or
not there is a $g \in \EP_2$ such that $g^{-1}yg=z$. }

\begin{proof}
A straightforward extension of Theorem 4.1 in \cite{matucci3} yields that $y \sim_{\EP_2}
z$ if and only if there exists $v_i \in T_{C_i}$ such that $v_i(t+1)=v_i(t)+1$, for $i=0,1$
and they satisfy equation (\ref{eq:ODP-mather}). Since $v_0$ needs to be equal to
$y^{-\infty}v_1 z^{\infty}$, our problem is reduced to deciding whether or not there is
$v_1 \in T_{C_1}$ solving these equations:
\begin{equation}\label{eq:ODP-equations-1}
\numberwithin{equation}{section}
\begin{array}{cc}
v_1(t+1)=v_1(t)+1, & \forall t \in C_1 \\
y^{-\infty}v_1 z^{\infty}(t+1)=y^{-\infty}v_1 z^{\infty}(t)+1, & \forall t \in C_0.
\end{array}
\end{equation}
Recalling that $C_0$ is a circle of length $a$ and $C_1$ is a circle of length $b$, we
define $s_i:C_i \to C_i$ to be the rotation by $1$ in $C_i$, for $i=0,1$. The problem now
becomes this: we need to decide whether or not there exists a map $v_1 \in T_{C_1}$ such
that
\begin{equation}\label{eq:ODP-equations-2}
\numberwithin{equation}{section}
\begin{array}{c}
v_1s_1=s_1v_1 \\
y^{-\infty} v_1 z^{\infty}s_0= s_0y^{-\infty}v_1 z^\infty.
\end{array}
\end{equation}
If we relabel $y^\infty s_0y^{-\infty}:=y^\ast$ and $z^{\infty}s_0 z^{-\infty}:=z^\ast$,
equations (\ref{eq:ODP-equations-2}) become
\begin{equation}\label{eq:simultaenous-equation}
\begin{array}{c}
v_1^{-1}s_1v_1 = s_1 \\
v_1^{-1}y^\ast v_1 =z^\ast.
\end{array}
\end{equation}
Equations (\ref{eq:simultaenous-equation}) are an instance of $2$-CP which is solvable by
assumption.
\end{proof}

\subsection{Non-trivial extensions of $F$ with solvable conjugacy problem}

\theorem{\label{thm:ODP-solvable}
If Conjecture~\ref{conj-T} is true for $k=2$, then $\Aut(F)$ and $\Aut_+(F)$ are orbit
decidable (as subgroups of $\Aut (F)$). In particular, assuming that such conjecture is
true, every group $G$ in an algorithmic short exact sequence
 $$
1 \longrightarrow F \overset{\alpha}{\longrightarrow} G \overset{\beta}{\longrightarrow} H \longrightarrow 1,
 $$
where $F=\PL_2(I)$, $H$ is a torsion-free hyperbolic group, and the action subgroup $A_G$
is either $\Aut(F)$ or $\Aut_+(F)$, has solvable conjugacy problem.
}

\begin{proof}
An application of Remark \ref{thm:remark-fixed-points-coincide} and Lemmas
\ref{thm:ODP-with-fixed-points} and \ref{thm:ODP-without-fixed} implies the solvability of
orbit decidability for the groups $\Aut(F)$ and $\Aut_+(F)$. We verify the requirements of
Theorem \ref{thm:bomave-extensions}. By Theorem \ref{thm:TCP-solvable}, condition (1) is
satisfied. It is well known (see, for example, Proposition 4.11(b) \cite{bomave2}) that if
$H$ is a free group or a torsion-free hyperbolic group, conditions (2) and (3) from Theorem
\ref{thm:bomave-extensions} are satisfied. By Theorem \ref{thm:ODP-solvable} we know that
the action subgroup is orbit decidable, then Theorem \ref{thm:bomave-extensions} implies
that $G$ has solvable conjugacy problem.
\end{proof}

\section{
\label{ssec:R-infty} Property $R_\infty$ in Thompson groups $F$ and $T$ }

In this section we show that Thompson groups $F$ and $T$ both have property $R_\infty$. We
recall the definition of property $R_\infty$, for the reader's convenience.

\definition{A group $G$ has property $R_\infty$ if for any $\varphi \in \Aut(G)$,
there exists a sequence $\{z_i\}_{i \in \mathbb{N}}$ of pairwise distinct elements which
are pairwise not $\varphi$-twisted conjugate. See also Section \ref{sec:intro}.}

We know that an automorphism $\varphi$ of $F$ is obtained by conjugation in $F$ by an
element $\tau\in\WTEP_2$. Moreover, we have seen in Subsection \ref{sec:restatement-TCP}
that two elements $y,z\in F$ are $\varphi$-twisted conjugate if and only if the two
elements $y\tau$ and $z\tau$ (now elements of $\WTEP_2$) are conjugate by an element of
$F$. Therefore, to prove that $F$ has property $R_\infty$ it is enough to show that, given
$\tau\in\WTEP_2$, there exists a family of elements $z_i\in F$, for all
$i=1,2,\ldots,n,\ldots$ such that they are pairwise not $\varphi$-twisted conjugate, i.e.,
$z_i\tau$ and $z_j\tau$ are not conjugate by an element of $F$.

Assume first that $\tau\in \EP_2$. If two elements are conjugate by an element of $F$ then
their fixed point sets match each other. So to prove that $z_i\tau$ and $z_j\tau$ are not
conjugate, it would be enough to construct the $z_i \in F$ in such a way that $z_i\tau$
has, say, a fixed point set with $i$ connected components so that the fixed point sets for
all the $z_i\tau$ would be different and the elements cannot be conjugate.

We observe that the fixed point set of $z_i\tau$ contains exactly the points
$t\in\mathbb{R}$ such that $z_i(t)=\tau^{-1}(t)$. Thus, it is enough to construct a map
$z_i\in F$ such that it has exactly $i$ disjoint intervals where $z_i(t)=\tau^{-1}(t)$,
thus producing $i$ connected components for $\Fix(z_i\tau)$. A reader familiar with $F$
should be able to construct easily such family $z_i$.

The proof above does not work if $\tau$ is orientation reversing. But it can be modified to
solve this case too. Assume now that $\tau=\sigma\mathcal{R}$ with $\sigma\in \EP_2$.
Construct the elements $z_i\in F$ similarly to the orientation preserving case using
$\sigma$, but in such a way that the fixed point set for $z_i\sigma$ is symmetric with
respect to the origin. More precisely, we can ensure that $\Fix(z_i\sigma)$ has $2i+1$
connected components given by $\{0\}$, $i$ connected components inside $\R_+$ and the
opposite of these components in $\R_-$. Moreover, we can ensure that $z_i \sigma > 0$ if
and only if $t>0$. Observe that by this symmetry, the map $\mathcal{R}z_i\sigma\mathcal{R}$
has the exact same fixed points as $z_i\sigma$ and so
$\Fix((z_i\sigma\mathcal{R})^2)=\Fix((z_i\sigma)^2)$.

Using this family $z_i$, we see that if $z_i\tau$ and $z_j\tau$ were conjugate via an
element of $F$, then $(z_i\sigma\mathcal{R})^2$ and $(z_j\sigma\mathcal{R})^2$ would also
be, and these have a different number of connected components in their fixed-point sets, by
construction, yielding a contradiction.

The argument above shows that we can recover property $R_\infty$ for $F$, giving a new
proof of the following result.

\theoremname{Bleak-Fel'shtyn-Gon\c{c}alves, \cite{bleakfelshgonc1}}{\label{thm:R-infty-F}
Thompson's group $F$ has property $R_\infty$. }

\remark{We notice that very recently Koban and Wong \cite{kowon} have shown that the group
$F \rtimes \mathbb{Z}_2$ has property $R_\infty$. }

Since we have a characterization for $\Aut(T)$ also in terms of conjugation by
piecewise-linear maps, the method described above to prove property $R_\infty$ for $F$ can
be used for $T$ as well.

\medskip
\noindent \textbf{Theorem \ref{thm:R-infty-T}.} \emph{Thompson's group $T$ has property
$R_\infty$.}

\begin{proof}
By Theorem 1 in \cite{brin5}, the group $\Aut(T)$ can be realized by inner automorphisms
and by conjugations by $\mathcal{R}$, the map which reverses the orientation.

The process will consist on constructing maps with different fixed-point sets. Consider a
piecewise-linear map on $[0,1]$ whose only fixed points are 0, $\frac{1}{2}$ and 1, and
also such that the graph is symmetric respect to the point $[\frac{1}{2},\frac{1}{2}]$.
Identify the endpoints to obtain a map on $S^1$ and hence an element of $T$. Call this map
$h_1$ and consider its lift $\widetilde{h}_1 \in \PL_2(\R)$. From the way we have
constructed $h_1$, we see that $\widetilde{h}_1$ is symmetric respect
$[\frac{1}{2},\frac{1}{2}]$ inside the square $[0,1]^2$, and so $\widetilde{h}_1$ is
invariant under $\mathcal{R}$, i.e.,
$\mathcal{R}\widetilde{h}_1\mathcal{R}=\widetilde{h}_1$ inside $\PL_2(\R)$. Therefore
$\mathcal{R}h_1 \mathcal{R} = h_1$ in $T$.

Now define inductively the map $h_i$ by subdividing the interval $[0,1]$ in its two halves
and in each half define a scaled-down version of $\widetilde{h}_{i-1}$, by a factor of 2.
Observe that if $i\neq j$, then $h_i$ and $h_j$ have different number of fixed points. For
a fixed $\varepsilon \in \{0,1\}$, if $h_i\mathcal{R}^{\varepsilon}$ and
$h_j\mathcal{R}^{\varepsilon}$ were conjugate in $T$, then
$(h_i\mathcal{R}^{\varepsilon})^2$ and $(h_j\mathcal{R}^{\varepsilon})^2$ are also
conjugate in $T$. We notice that $(h_i\mathcal{R})^2=h_i^2$ and that $h_i^2$ and $h_j^2$
have different number of fixed points, so they cannot be conjugate.
\end{proof}

\section{Generalizations and some questions
\label{sec:generaltions-of-results}}

In this section we make a series of observations about the extent to which the material of
this paper generalizes and describe some natural related questions.

\subsection{Extensions of the Bieri-Thompson-Stein-Strebel groups $\PL_{S,G}(I)$} It
seems likely that the theory developed in this paper can be generalized to a certain extent
to the Bieri-Thompson-Stein-Strebel groups $\PL_{S,G}(I)$, with the computational
requirements described in \cite{matucci5}.

We recall that $\PL_{S,G}(I)$ is the group of piecewise-linear homeomorphisms of the unit
interval $I$ with finitely many breakpoints occurring inside $S \leqslant \mathbb{R}$, an
additive subgroup of $\mathbb{R}$ containing $1$, and such that the breakpoints lie in $G
\leqslant U(S)$, where $U(S)=\{g \in \mathbb{R}^* \mid gS=S \text{ and } g>0 \}$.

Since our results rely on straightforward generalizations of those in \cite{matucci5} and
\cite{matucci3}, to generalize our algorithms to the groups $\PL_{S,G}(I)$ we need to
observe a number of things:
\begin{enumerate}
\item We define the analogues $\PL_{S,G}(\mathbb{R}),\WTEP_{S,G}, \EP_{S,G}$ and observe
    that the existence of periodicity boxes,  the construction of conjugators and moving
    fixed points (Subsections \ref{ssec:periodicity-boxes}, \ref{ssec:periodicity-boxes}
    and \ref{ssec:fixed-points}) generalize immediately via the results in
    \cite{matucci5} (which are proved in $\PL_{S,G}(I)$).
\item To reduce the number of possible ``initial slopes'' we need to generalize
    Subsection \ref{sec:rescaling-the-circle}. We can do this since the material in
    \cite{matucci3} can be generalized to $\PL_{S,G}(I)$. The second observation that is
    needed to reduce slopes is the one used in the proof of Theorem
    \ref{thm:TCP-solvable}, where we multiply a candidate conjugator $g$ by a power of
    $y^2$. This shows that we need to build candidate conjugators only for slopes in
    $[(y^2)'(p^+),1]$ and, by Lemma 5.4 in \cite{matucci5}, we can show that the sets of
    slopes is discrete in $\mathbb{R}_+$, thereby giving us only finitely many slopes
    inside $[(y^2)'(p^+),1]$. Hence, this part generalizes too.
\item Brin's Theorem \ref{thm:brin-thm} has a non-trivial generalization in a result of
    Brin and Guzman \cite{bringuzman} which describes certain classes of automorphisms of
    the groups $\PL_{\mathbb{Z}[\frac{1}{n}],\langle n \rangle}(I)$. There exist elements
    in the automorphism group $\Aut(\PL_{\mathbb{Z}[\frac{1}{n}],\langle n \rangle}(I))$
    which are represented by conjugation via elements that are not in $\WTEP_n$ (and that
    are called ``exotic''). Therefore, we can only generalize results of the current
    paper by restricting the action subgroup being used. Instead of studying the full
    automorphism group $\Aut(\PL_{S,G}(I))$, we can restrict to study conjugations by
    element of $\WTEP_{S,G}$ so that we can adapt our results in a straightforward
    manner.
\end{enumerate}

\remark{It should be noted that the tools of this paper are not generally sufficient to
solve either the twisted conjugacy problem or the orbit decidability problem in any group
$\PL_{S,G}(I)$ generalizing Thompson's group $F$ (for example, in generalized Thompson's
groups $F(n)$). This is because the full automorphism group may contain conjugations via
not piecewise-linear maps.}

It is however possible to give suitable reformulations of Theorems ~\ref{thm:TCP-solvable},
~\ref{thm:ODP-solvable} and
~\ref{thm:CP-extension-unsolvable} in the setting of actions whose acting group is realized
by conjugations by an element of $\WTEP_{S,G}$. The restatement of Theorem
~\ref{thm:ODP-solvable} will need to assume that the
$2$-simultaneous conjugacy problem is solvable for the groups $T_{S,G}$ and this is also
work-in-progress as mentioned in Section ~\ref{sec:ODP-solvable}.

Since the techniques used to study the twisted conjugacy problem for $F$ arise from those
used in \cite{matucci5} to study the simultaneous conjugacy problem for $F$, it is natural
to ask the following question:

\question{Is the $k$-simultaneous twisted conjugacy problem solvable for $F$? More
precisely, is it decidable to determine whether or not, given $\varphi \in \Aut(F)$ and
$y_1,\ldots,y_k,z_1,\ldots,z_k \in F$, there exists a $g \in F$ such that $z_i=g^{-1}y_i
\varphi(g)$?}

\subsection{Extensions of Thompson's group $T$}
As observed at the beginning of the proof of Theorem \ref{thm:R-infty-T}, if $\varphi \in
\Aut(T)$, then there exists an $\varepsilon \in \{0,1\}$ such that
$\varphi(\lambda)=\mathcal{R}^{\varepsilon} \tau^{-1} \alpha \tau
\mathcal{R}^{\varepsilon}$, for all $\alpha \in T$. Arguing as in Subsection
\ref{sec:restatement-TCP}, equation (\ref{eq:TCP-equation}) can be rewritten as
\begin{equation}
\label{eq:TCP-in-T-again} \numberwithin{equation}{section} g^{-1}(y\mathcal{R}^\varepsilon) g =
z\mathcal{R}^\varepsilon
\end{equation}
for $y,z,g \in T$ and $\varepsilon \in \{0,1\}$. To attack equation
(\ref{eq:TCP-in-T-again}), we can start by squaring it and initially reduce ourselves to
solve the equation
\begin{equation}
\label{eq:squared-TCP-in-T-again} \numberwithin{equation}{section} g^{-1}(y\mathcal{R}^\varepsilon)^2 g =
(z\mathcal{R}^\varepsilon)^2.
\end{equation}
The advantage of working with equation (\ref{eq:squared-TCP-in-T-again}) is that
$(y\mathcal{R}^\varepsilon)^2, (z\mathcal{R}^\varepsilon)^2 \in T$.

The conjugacy problem in $T$ is solvable by the work of Belk and the second author in
\cite{matucci9} and thus we can list all the conjugators in $T$ between
$(y\mathcal{R}^\varepsilon)^2$ and $(z\mathcal{R}^\varepsilon)^2 $. However, there might be
infinitely many of them and there is no obvious way to detect which of them will also be
conjugators between $y\mathcal{R}^\varepsilon$ and $z\mathcal{R}^\varepsilon$.

We cannot use the techniques of the current paper, since there is no uniqueness given by an
the ``initial slope'' of elements of $T$ (although something similar may be feasible, as it
was done in Chapter 7 in \cite{matuccithesis} to study centralizers in $T$). We are thus
led to ask:

\question{Is the twisted conjugacy problem solvable in Thompson's group $T$?}

To conclude, we mention that the orbit decidability problem for $T$ is solvable for $\Aut(T)$
and $\Aut_+(T)$.

\lemma{Let $T$ be Thompson's group $\PL_2(S^1)$. Then $\Aut(T)$ and $\Aut_+(T)$ are orbit
decidable.}

\begin{proof}
We need to decide whether or not, given $y,z \in T$, there exists an element $g \in T$ such
that at least one of the two equalities
\begin{equation}
\label{eq:ODP-in-T} \numberwithin{equation}{section} z=g^{-1}yg \qquad \text{or} \qquad
z=g^{-1}(\mathcal{R}y\mathcal{R})g
\end{equation}
holds. This amounts to study two distinct conjugacy problems for elements of $T$, each of
which is solvable by the work \cite{matucci9}.
\end{proof}

\bibliographystyle{plain}
\bibliography{go4}

\end{document}